\documentclass{amsart}
\usepackage{amsmath, amsthm}
\usepackage{amssymb, graphicx}

\theoremstyle{plain}
\newtheorem{prop}{Proposition}[section]
\newtheorem*{propA}{Proposition A}
\newtheorem*{propB}{Proposition B}
\newtheorem{coro}[prop]{Corollary}
\newtheorem{lemm}[prop]{Lemma}
\newtheorem{thrm}[prop]{Theorem}
\newtheorem*{thrmA}{Theorem A}
\newtheorem*{thrmB}{Theorem B}

\theoremstyle{definition}
\newtheorem{defi}[prop]{Definition}

\newtheorem{exam}[prop]{Example}
\newtheorem{rema}[prop]{Remark}

\numberwithin{equation}{section}

\def\System#1{\pmb{\mathsf{#1}}}
\renewcommand\a{\alpha}%
\newcommand\Ack{\mathrm{Ack}}%Ackermann
\newcommand\BB[1]{B_{#1}^{\raise1.5pt\hbox{$\scriptscriptstyle+$}}}
\newcommand\BBi{\BB\infty}
\renewcommand\b{\beta}   
\newcommand\bb{b}   
\newcommand\can[1]{\vert\!\vert#1\vert\!\vert}%canonical_length
\newcommand\card{\mathrm{card}}%cardinal
\newcommand\CC[3]{\langle#3\rangle_{_{\!#1,#2}}}%compositiontOfSpecialBraids
\newcommand\cl[1]{[#1]}%parity
\newcommand\Const[1]{c_{#1}}
\newcommand\D{\Delta}
\newcommand\dd{d}
\newcommand\ddd[1]{\delta_{#1}}
\newcommand\dddd[1]{\theta_{#1}}
\newcommand\Div{\mathrm{Div}}
\newcommand\emin{e^{\hbox{\smaller\smaller\smaller\smaller min}}}%minimalExponent
\newcommand\ee{e}%exponent
\newcommand\eps{\varepsilon}%ordinal
\newcommand\eqCNF{=_{\scriptscriptstyle\mathrm{CNF}}}%CantorNormalForm

\newcommand\expo[1]{d(#1)}%exponent_of_Delta
\newcommand\g{\gamma}
\renewcommand\ge{\geqslant}
\newcommand\ff{f}%complexity_function
\newcommand\flip[1]{\phi_{#1}}%flip
\newcommand\flipt[1]{\widetilde{\phi}_{#1}}%special_flip
\newcommand\For{\Phi}%formula
\newcommand\Forr{\Psi}%formula
\newcommand\G{\mathcal{G}}%game
\renewcommand\gg{g}%function
\newcommand\Gsp{\mathcal{G}^{\raise1.5pt\hbox{$\scriptscriptstyle sp$}}}%game
\newcommand\hh{h}%function

\newcommand\IDEXP{\System{I\Delta}_0{+}\boldsymbol{e\!x\!p}}
\newcommand\ie{{\it i.e.}}
\newcommand\ii{i}
\newcommand\Ind{\boldsymbol{I\!n\!d}}
\newcommand\inv{^{-1}}
\newcommand\IS[1]{\System{I\Sigma}_{#1}}
\newcommand\jj{j}
\newcommand\kk{k}
\renewcommand\l{\lambda}
\renewcommand\le{\leqslant}
\renewcommand\lg[1]{\vert#1\vert}%word_length
\newcommand\logb[1]{\ell\hspace{-0.3ex}o\hspace{-0.2ex}g{#1}}%word_length
\newcommand\lSL{<^{\scriptscriptstyle\mathtt{ShortLex}}}
\newcommand\mm{m}
\newcommand\Nat{\mathbb{N}}
\newcommand\nn{n}
\newcommand\nnd{{n-2}}
\newcommand\nno{{n-1}}

\newcommand\om{\omega}
\newcommand\op{\cdot}
\newcommand\OP[1]{\{#1\}}
\newcommand\OPsp[1]{\{#1\}^{\raise1pt\hbox{$\scriptscriptstyle\! sp$}}}
\newcommand\ord{\mathrm{ord}}
\newcommand\ordsp{\mathrm{ord}^{\raise1pt\hbox{$\scriptscriptstyle sp$}}}
\newcommand\PA{\System{PA}}
\newcommand\pp{p}
\newcommand\PRA{\System{PRA}}
\newcommand\qq{q}
\newcommand\resp{{\it resp.}~}
\newcommand\rr{r}

\newcommand\seq[2]{c_{#1,#2}}%
\newcommand\Seq[1]{\Sigma_{#1}}
\newcommand\Seqi[1]{\sigma_1^{(#1)}}
\newcommand\Seqt[2]{\widetilde\Sigma_{#1,#2}}
\newcommand\ShortLex{\mathtt{ShortLex}}
\renewcommand\ss[1]{\sigma_{\!#1}^{\vrule height5pt width0pt}}%braidSigma 
\renewcommand\SS{S}%set
\newcommand\sss[2]{\sigma_{\!#1}^{#2}}%braidSigma 

\renewcommand\tt{t}
\newcommand\TT{T}
\newcommand\TTsp{T^{\raise1.5pt\hbox{$\scriptscriptstyle sp$}}}
\newcommand\tti{s}
\newcommand\ttii{t}
\newcommand\UU{U}%function
\newcommand\UUsp{U^{\raise1.5pt\hbox{$\scriptscriptstyle sp$}}}%function
\newcommand\WO[1]{\boldsymbol{W\!O}_{\!#1}}%combinatorial_wf
\newcommand\ww{w}%aBraidWord
\newcommand\xx{x}%
\newcommand\yy{y}%aBraid
\newcommand\ZF{\System{ZF}}
\newcommand\zz{z}%

\begin{document}

\author{Lorenzo CARLUCCI}
\address{Universit\`a di Roma ``La Sapienza'',
Department of Computer Science,
Via Salaria 113, 00198 Roma, Italy\\and\\
Scuola Normale Superiore di Pisa,
Classe di Lettere, Piazza dei Cavalieri, 56126 Pisa, Italy}
\email{carlucci@di.uniroma1.it}

\author{Patrick DEHORNOY}
\address{Laboratoire de Math\'ematiques Nicolas Oresme,
UMR 6139 CNRS,  Universit\'e de Caen BP 5186, 14032 Caen,
France}
\email{dehornoy@math.unicaen.fr}
\urladdr{//www.math.unicaen.fr/\!\hbox{$\sim$}dehornoy}

\author{Andreas WEIERMANN}
\address{Vakgroep Zuivere Wiskunde en Computeralgebra
Ghent University, Krijgslaan 281, Gebouw S22,
B9000 Gent, Belgium}
\email{weierman@cage.ugent.be}

\title{Unprovability results involving braids}


\thanks{L.\,C. was partially supported  by Telecom Italia ``Progetto Italia" grant; L.\,C.  and A.\,W. were 
supported by NWO grant number  613080000.}


\begin{abstract}
We construct long sequences of braids that are descending
with respect to the standard order of braids (``Dehornoy
order''), and we deduce that, contrary to all usual algebraic properties of braids,
certain simple combinatorial statements involving the braid order are true, but not
provable in the subsystems~$\IS1$ or~$\IS2$ of the
standard Peano system. 
\end{abstract}

\maketitle

It has been known for decades that there exist strong limitations about the
sentences possibly provable from the axioms of a given formal system,
starting with G\"odel's famous theorems implying that certain arithmetic
sentences cannot be proved from the axioms of the first-order Peano system.
However, the so-called G\"odel sentences have a strong logical flavour and
they remain quite remote from the sentences usually considered by
mainstream mathematicians. It is therefore natural to look for further
sentences that are true but unprovable from the Peano axioms, or from the
axioms of other formal systems, and, at the same time, involve objects and
properties that are both simple and natural. The main results so far in this direction involve finite combinatorics, Ramsey Theory and the theory of well-quasi-orders. See~\cite{Bov, Wei2} for a comprehensive bibliography.

On the other hand, Artin's braid groups are algebraic structures which play  a
central role in many areas of mathematics and theoretical
physics~\cite{Bir, KaT}. It has been known since 1992 that, for each~$\nn \ge
2$, the group~$B_\nn$ of $\nn$-strand braids is equipped with a canonical
left-invariant ordering~\cite{Dfb}, and one of the most remarkable properties
of this ordering is the result, due to R.\,Laver~\cite{Lve}, that its restriction to
the submonoid~$\BB\nn$ of~$B_\nn$ consisting of the so-called Garside
positive braids is a well-order, \ie, every nonempty subset of~$\BB\nn$ has a least
element. It was proved by S.\,Burckel in~\cite{Bur} that the
order-type of this restriction is the ordinal~$\om^{\om^{\nn-2}}$, hence it is
rather large in the hierarchy of well-orders. It follows that, although the existence of infinite
descending sequences in~$\BB\nn$ is forbidden by the well-order property, there may exist {\it
long} finite descending sequences. 

What we do in this paper is to investigate the existence of such long descending
sequences in~$(\BB\nn, <)$ from the viewpoint of provability in~$\IS1$ and~$\IS2$,  the
subsystems of the Peano system in which the induction scheme restricted to $\Sigma_1$
and $\Sigma_2$ sentences respectively, where a sentence is $\Sigma_\kk$
if it is of the form $\exists x_1\forall x_2\exists x_3\,...\, Q \xx_\kk(\Phi)$,
with $\kk$ quantifiers and $\Phi$ containing bounded quantifiers only---see Appendix for complete
definitions; more generally, the few notions from logic needed for the paper are recalled there. We establish
two types of unprovability results, that we now state in the context of~$\BB3$,
\ie, of $3$-strand braids. First, we introduce particular long descending sequences of braids,
called $\G_3$-sequences, by a simple recursive process. Then we prove

\begin{propA}
For each initial braid~$\bb$ in~$\BB3$, the $\G_3$-sequence from~$\bb$
is finite.
\end{propA}

\begin{thrmA}
Proposition~A is an arithmetic statement\/\footnote{By an arithmetic statement we
mean a first-order sentence in the language of Peano Arithmetic. As it stands, Proposition~A
involves braids, and therefore it is not an arithmetic statement; what we mean
is that Proposition~A can be encoded into an arithmetic statement in a way
whose correctness can be established using the axioms of~$\IS1$.} that is not provable from the axioms
of~$\IS1$.
\end{thrmA}

By contrast, it should be emphasized that much of the usual properties of
braids, in particular all known algebraic properties, can, when properly
encoded, be proved from the axioms of~$\IS1$---as do most of the usual
mathematical results that are formalizable in that system.

The second family of results involves general descending sequences of
braids, and not only those called $\G_3$-sequences in Proposition~A. For
each function~$\ff$ of~$\Nat$ to~$\Nat$, we introduce a certain combinatorial
principle~$\WO\ff$ that, roughly speaking, says that each descending sequence in~$\BB3$ in
which the Garside complexity of the $\kk$th braid entry remains below~$\ff(\kk)$ has a
bounded length. We establish

\begin{propB}
For each function~$\ff$, the principle~$\WO\ff$ is true.
\end{propB}
But, denoting by~$\Ack$ the standard Ackermann function---see
Appendix---and by~$\Ack_\rr$ the level~$\rr$ approximation to~$\Ack$, and using
$\ff\inv$ for the functional inverse of~$\ff$, we prove

\begin{thrmB}
$(i)$ For~$\rr \ge 0$, let $\ff_\rr$ be defined by $\ff_\rr(\xx) =
\lfloor\!\!\sqrt[\Ack_\rr\inv(\xx)]{\xx}\rfloor$, and $\ff_\om$ be defined by 
$\ff_\om(\xx) =  \lfloor\!\!\sqrt[\Ack\inv(\xx)]{\xx}\rfloor$. 

$(i)$ For each~$\rr$, the principle~$\WO{\ff_\rr}$ is provable from the axioms of~$\IS1$.

$(ii)$ The principle~$\WO{\ff_\om}$ is not provable from the axioms of~$\IS1$. 
\end{thrmB}

The functions involved in Theorem~B all are of the form $\xx\mapsto
\!\!\sqrt[f(\xx)]\xx$ where $f$ is a very slowly increasing function.
Analogous to the results of~\cite{Wei2, Wei3, Wei4}, Theorem~B is a typical
example of a so-called phase transition phenomenon, in which a seemingly
small change of the parameters causes a jump from provability to
unprovability, here with respect to~$\IS1$.

In some sense, the above results about $3$-strand braids, as well as their
extensions involving arbitrary braids, are not surprising. As $3$-strand braids (\resp general
braids) are equipped with a well-ordering of length~$\om^\om$ (\resp $\om^{\om^\om}$),  a
connection with the system~$\IS1$ (\resp $\IS2$) can even be expected, because of the
well-known connection of the latter ordinal with that logical system---{\it cf.} for instance Simpson's analysis 
of the Hilbert and Robson basis theorems in~\cite{Sim}. The results we establish are reminiscent of
analogous results established in the language of ordinals and trees. For instance, our  $\G_3$-sequences are
direct cousins of the Goodstein sequences and  the Hydra battles~\cite{KiP} as well as of the more recent Worm
Principle~\cite{Bek, Lee}. However, our results are not just artificial translations of existing properties into the
language of braids. The braid order is arguably a quite natural object, and all arguments  developed in this
paper rely on the specific properties of braids and their order, and not on an automated translation into
another context. Typically, Propositions~A and~B directly follow from the very definition of the braid
ordering and its well-foundedness, while Theorems~A and~B rely on some non-trivial analysis of the braid
order on~$\BB3$ and its connection with Garside's theory. The reason that makes the current results
essentially nontrivial is that, although the well-order on positive braids is just a copy of the well-order on
ordinals---according to the general uniqueness theorem of well-orders of a given length---the actual order
isomorphism between braids and ordinals is not simple. This explains in particular why relatively
sophisticated braid arguments are needed here. At the very least, one of the interests in the current approach
is that it led to interesting braid questions that could be solved only at the expense of developing new tools,
such as the counting results of Section~\ref{SS:Counting} or the decomposition results of
Section~\ref{S:Special}.

The paper is organized as follows. Section~\ref{S:Context} contains a brief
introduction to braid groups, their ordering, and to the so-called
$\flip{}$-normal form of $3$-strand braids, all needed to state the subsequent results. In
Section~\ref{S:Sequence3}, we describe the
$\G_3$-sequences, and establish Proposition~A and Theorem~A. In
Section~\ref{S:Transition}, we introduce the combinatorial
principle~$\WO\ff$, and establish Proposition~B and Theorem~B. Finally, in
Section~\ref{S:Extension}, we show how to extend Proposition~A and Theorem~A into similar
results involving general braids and $\IS2$-provability. We also raise a few questions and point
to further research. Finally, we provide in an appendix the needed basic definitions from
logic, about ordinals and about basic subsystems of Peano arithmetic.

\section{The general braid context}
\label{S:Context}

We briefly recall definitions for the braid group~$B_\nn$,
the braid monoid~$\BB\nn$, and the canonical braid order that will be the
central object of investigation in the subsequent sections. The main point is
the connection between the braid order on~$\BB3$ and the so-called
$\flip{}$-normal form.

\subsection{Braid groups}
\label{SS:Braids}

For $\nn\ge2$, the $\nn$-strand braid group~$B_n$ is the 
group of isotopy classes of geometric $\nn$-strand
braids~\cite{Bir, KaT}. For our current purpose, it is sufficient to
know that
$B_\nn$ is the group with presentation
\begin{equation}
\label{E:Pres}
 \langle \ss1, ... , \ss{\nn-1} \, ; \,
 \ss i \ss j = \ss j \ss i \; 
 \mbox{for $\vert i - j\vert \ge 2$,}\; 
 \ss i \ss j \ss i = \ss j \ss i \ss j 
 \; \mbox{for $\vert i - j\vert = 1$}
 \rangle,
\end{equation}
so that every element of~$B_\nn$, called an $\nn$-braid in
the sequel, is an equivalence class of words  on the
letters~$\sss1{\pm1}$, ..., $\sss{\nn-1}{\pm1}$ with respect to the
congruence generated by the relations of~\eqref{E:Pres}. 

\begin{figure} [htb]
\begin{picture}(82,15)(0,-3)
\put(0,0){\includegraphics{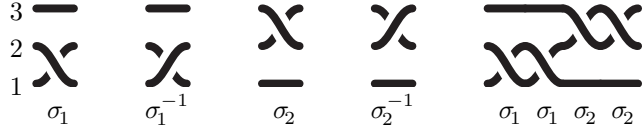}}
\put(-3,-1){$1$}
\put(-3,4){$2$}
\put(-3,9){$3$}
\put(2,-4){$\ss1$}
\put(15,-4){$\sss1{-1}$}
\put(32,-4){$\ss2$}
\put(45,-4){$\sss2{-1}$}
\put(62,-4){$\ss1$}
\put(67,-4){$\ss1$}
\put(72,-4){$\ss2$}
\put(77,-4){$\ss2$}
\end{picture}
\caption{\smaller\sf Three strand braid diagrams 
associated with~$\ss1$, $\sss1{-1}$, $\ss2$, $\sss2{-1}$, and
$\sss12\sss22$: positions are numbered $1$ to $3$ from
bottom, $\ss\ii$ (\resp $\sss\ii{-1}$) denotes the half-twist
where the strand at position~$\ii+1$ crosses over (\resp
under) the strand at position~$\ii$; the diagram associated
with a product~$\ww\ww'$ is the concatenation of
the diagrams associated with~$\ww$ and~$\ww'$.}
\label{F:Sigma}
\end{figure}

The connection with geometry is as follows. Associate with every $\nn$-strand 
braid word~$\ww$ a braid diagram by concatenating the elementary diagrams
of Figure~\ref{F:Sigma} corresponding to the successive letters of~$\ww$. Such
a diagram can be seen as a plane projection of a three-dimensional figure
consisting on $\nn$~disjoint curves. Then, the relations of~\eqref{E:Pres}
are a translation of ambient isotopy, \ie, of continuously moving
the curves without moving their ends and without allowing
them to intersect. It is easy to check that the relations
of~\eqref{E:Pres} correspond to such isotopies; the converse
implication, \ie, the fact that the projections of isotopic
three-dimensional geometric braids always can be encoded
in words connected by~\eqref{E:Pres}, was proved by E.~Artin
in~\cite{Art}. 

For each~$\nn$, the identity mapping on~$\{\ss1,...,
\ss{\nn-1}\}$ induces an embedding of~$B_\nn$
into~$B_{\nn+1}$, and we shall henceforth identify~$B_\nn$
with its image in~$B_{\nn+1}$---so that there is no need to
distinguish the generator~$\ss i$ of~$B_\nn$ from that
of~$B_{\nn+1}$. Geometrically, this corresponds to freely
adjoining  additional unbraided strands on the top of diagrams. 
In this way, the groups~$B_\nn$ arrange into a direct system; 
its direct limit is denoted~$B_\infty$, and it is the group
generated by an infinite sequence of generators $\ss1,
\ss2, ...$ subject to the relations of~\eqref{E:Pres}.

Positive braids are defined to be those braids that admit an
expression involving no letter~$\sss\ii{-1}$. For $\nn\le\infty$,
positive $\nn$-braids make a submonoid~$\BB\nn$
of~$B_\nn$ that is known to admit, as a monoid, the presentation~\eqref{E:Pres}.

\subsection{The flip normal form}

As several braid diagrams may be associated with a braid, \ie, several
braid words may represent the same braid, it it is often convenient to
distinguish a particular representative in each equivalence class of braid
words. Several  classical solutions exist, in particular the so-called Garside, or
greedy, normal form~\cite[Chap.~9]{Eps}. In this paper, we shall briefly use
the greedy normal form in Section~\ref{S:Complexity} below, but our main
tool will be another normal form called the flip normal, or $\flip{}$-normal,
form, which was introduced by S.\,Burckel in~\cite{Bur} and further
investigated in~\cite{Dho}. For our current purpose, it is enough to consider
the case of~$\BB3$, which is very simple.

A positive $3$-strand braid is unambiguously specified by a
word on the alphabet~$\{\ss1, \ss2\}$, \ie, a finite sequence of~$\ss1$'s
and~$\ss2$'s. Every such word~$\ww$, simply called {\it a $3$-braid word}
in the  sequel, can be expressed as
$\sss{\cl\pp}{\ee_\pp} \,...\, \sss1{\ee_3} \sss2{\ee_2}  \sss1{\ee_1}$, where
$\cl\pp$ means~$1$ for $\pp$~odd, and~$2$ for $\pp$~even. Among all
such expressions of~$\ww$, there exists a unique one in which $\pp$ has the
minimal value; it will be called the {\it block decomposition} of~$\ww$, and
the subwords~$\sss\kk{\ee_\kk}$---as well as the associated fragments in the
corresponding braid diagram---will be called the {\it blocks} of~$\ww$.

\begin{defi}
\label{D:Normal3}
Let $\emin_1 = 0$, $\emin_2 = 1$, and $\emin_\kk = 2$ for $\kk \ge 3$. A
$3$-braid word is said to be {\it $\flip{}$-normal} if its block decomposition
$\sss{\cl\pp}{\ee_\pp} \, ...\, \sss2{\ee_2}\sss1{\ee_1}$ satisfies the
inequalities $\ee_\kk \ge \emin_\kk$ for~$\kk  < \pp$.
\end{defi}

For instance, the word~$\ss1\ss2\ss1$ is $\flip{}$-normal, but the 
word~$\ss2\ss1\ss2$, whose block decomposition is
$\sss21\sss11\sss21\sss10$, is not, since the condition $\ee_3 \ge \emin_3$
fails: $\ee_3$ is~$1$ here. Using the braid relation $\ss1\ss2\ss1 =
\ss2\ss1\ss2$, one easily establishes:

\begin{prop}
\cite{Bur}
\label{P:Normal3}
Every nontrivial braid in~$\BB3$ can be represented by a unique
$\flip{}$-normal word.
\end{prop}

The name ``$\flip{}$-normal'' refers to the flip automorphism~$\flip3$
of~$\BB3$ that exchanges~$\ss1$ and~$\ss2$: it is shown in~\cite{Dho} that
the $\flip{}$-normal expression of a braid~$\bb$ can be obtained by
considering the maximal power of~$\ss1$ that is a right divisor of~$\bb$ in
the monoid~$\BB3$, then the maximal power of~$\ss2$ that is a right divisor
of the quotient, and on so considering~$\ss1$ and~$\ss2$, \ie, $\ss1$
and~$\flip3\ss1$, alternatively. 

\begin{defi}
\label{D:ExpSeq3}
Let $\bb$ be a nontrivial braid in~$\BB3$, and let $\sss{\cl\pp}{\ee_\pp} \,
...\, \sss2{\ee_2}\sss1{\ee_1}$ be the $\flip{}$-normal expression of~$\bb$ given
by Proposition~\ref{P:Normal3}. Then the sequence $(\ee_\pp, ..., \ee_1)$ is
called the {\it exponent sequence} of~$\bb$, and the number~$\pp$ is called
its {\it breadth}. 
\end{defi}

For instance, let~$\D_3$ be the braid~$\ss1 \ss2 \ss1$---the so-called Garside's
fundamental $3$-braid of~Ê\cite{Gar}. We observed above that the word $\ss1 \ss2 \ss1$ is
$\flip{}$-normal, so the exponent sequence of~$\D_3$ is~$(1, 1, 1)$, and its
breadth is~$3$. Similarly, it is easily checked~\cite{Dho} that, for
each~$\kk$, the exponent sequence of~$\D_3^\kk$ is $(1, 2, ..., 2, 1, \kk)$,
with $2$ repeated $\kk-1$~times; thus the breadth of~$\D_3^\kk$
is~$\kk+2$.

\subsection{The braid order}

Braids are equipped with a canonical linear order. The latter can be simply defined
in terms of word representatives of a specific form---but it also admits a number of equivalent
definitions~\cite{Dgr}.

\begin{defi}
\label{D:Order}
If $\bb, \bb'$ are braids, we say that $\bb < \bb'$ holds if
the braid~$\bb\inv\bb'$ admits an expression by a braid word in which the
generator~$\ss\ii$ with higher index occurs only positively, \ie, $\ss\ii$
occurs, but~$\sss\ii{-1}$ does not\footnote{The current relation is
denoted~$<^{\!\scriptscriptstyle\phi}$ in~\cite{Dgr}: this version of the
braid order, which was first considered by S.\,Burckel in~\cite{Bur}, is
definitely more suitable when well-order properties are involved than the
symmetric version in which one takes into account the
generator~$\ss\ii$ with lowest index.}.
\end{defi}

For instance, $\ss2<\ss1\ss2$ holds, as the quotient
~$\sss2{-1}\ss1\ss2$ can also be expressed
as~$\ss1\ss2\sss1{-1}$, and, in the latter word, the main 
generator, which is~$\ss2$, occurs positively (one~$\ss2$), but
not negatively (no~$\sss2{-1}$). 

\begin{thrm}
\label{T:Order}
$(i)$ \cite{Dfb} The relation~$<$ is a linear ordering 
on~$B_\infty$ that is left compatible with multiplication;
for each~$\nn$, the set~$B_\nn$ is an open interval of~$(B_\infty,<)$ centered on~$1$.

$(ii)$ \cite{Lve, Bur} The restriction of~$<$ to~$\BB\infty$ is a well-ordering
of ordinal type~$\om^{\om^\om}$; for each~$\nn$, the restriction of~$<$ to~$B_\nn^+$ is the
initial segment~$[1, \ss\nn)$ of~$(\BB\infty,<)$, and its ordinal type is~$\om^{\om^{\nn-2}}$.
\end{thrm}

The braid order is nicely connected with the $\flip{}$-normal
form~\cite{Bur, Dho}. In the current paper, we shall be mostly dealing with~$\BB3$, and we can easily
describe, and even reprove, the connection in this simple case.

\begin{lemm}
\label{L:Order}
For $\pp \ge 0$, let $\ddd\pp$ be the braid represented by the length~$2\pp$ suffix of the left
infinite word $...\sss12\sss22\sss12\ss2$. Then, for $\pp \ge 1$, we have $\bb < \ddd\pp$ for each
braid~$\bb$ with breadth at most~$\pp + 1$, and $\ddd\pp \le \bb$ for each
braid~$\bb$ with breadth at least~$\pp + 2$.
\end{lemm}

\begin{proof}
Let $\bb$ be a braid in~$\BB3$ with breadth at most~$\pp+1$.
Then, by definition, we have $\bb = \sss{\cl{\pp+1}}{\ee_{\pp+1}} ...
\sss2{\ee_2} \sss1{\ee_1}$ for some nonnegative exponents $\ee_{\pp+1},
..., \ee_1$. An easy induction gives for each~$\pp \ge 0$ the equality
\begin{equation}
\label{E:Deltadelta}
\D_3^\pp = \ddd\pp \sss1\pp.
\end{equation}
So, for $\pp \ge 1$, we obtain 
\begin{equation}
\label{E:Delta}
\bb\inv \op \ddd\pp 
= \sss1{- \ee_1} \sss2{-\ee_2} \,...\, \sss{\cl{\pp+1}}{-\ee_{\pp+1}} \op
\D_3^\pp \, \sss1{-\pp}.
\end{equation}
The braid relations~\eqref{E:Pres} imply  $\ss\ii \op \D_3 =
\D_3 \op \flip3\ss\ii$ for $\ii = 1, 2$, where $\flip3$ is the automorphism
of~$\BB3$ that exchanges~$\ss1$ and~$\ss2$. This enables us to push the
factors~$\D_3$ of~\eqref{E:Delta} to the left, at the expense of
applying~$\flip3$. In this way, we deduce
\begin{equation*}
\bb\inv \op \ddd\pp 
= \sss1{- \ee_1} \op \D_3 \op \sss1{-\ee_2} \op \D_3 \op ... \op \D_3 \op 
\sss1{-\ee_{\pp+1}} \op \sss1{-\pp}.
\end{equation*}
whence, using~$\D_3 = \ss1\ss2\ss1$,
\begin{equation*}
\bb\inv \op \ddd\pp 
= \sss1{- \ee_1+1} \, \ss2 \, \sss1{-\ee_2+2} \, \ss2 \, ... \, \ss2 \, 
\sss1{-\ee_\pp + 2} \, \ss2 \, \sss1{-\ee_{\pp+1} + 1} \, \sss1{-\pp}.
\end{equation*}
The generator~$\ss2$ occurs $\pp$~times in the above expression, while  $\sss2{-1}$ does not
occur. Hence, by definition, $\bb <
\ddd\pp$ holds.

Assume now that $\bb$ has breadth $\pp + 2$ or more. Owing to
Proposition~\ref{P:Normal3}, we can write $\bb = \bb'
\sss{\cl{\pp+2}}{\ee_{\pp+2}} ... \sss2{\ee_2} \sss1{\ee_1}$ with
$\ee_{\pp+2} \ge 1$, $\ee_{\pp+1}, ..., \ee_3 \ge 2$, $\ee_2 \ge1$, and
$\ee_1 \ge 0$. We find
\begin{equation*}
\ddd\pp\inv \op \bb 
= \sss1\pp \, \D_3^{-\pp} \op \bb' \,
\sss{\cl{\pp+2}}{\ee_{\pp+2}} \, ... \, \sss2{\ee_2} \,\sss1{\ee_1}.
\end{equation*}
Pushing the factors~$\D_3\inv$ to the right, and using $\D_3\inv =
\sss2{-1} \sss1{-1} \sss2{-1}$, we deduce
\begin{align}
\notag
\ddd\pp\inv \op \bb
&= \sss1\pp \op \flip3\bb' \op 
\sss2{\ee_{\pp+2}} \op \D_3\inv \op \sss2{\ee_{\pp+1}} \op \D_3\inv \op ...
\op \D_3\inv \op  \sss2{\ee_2} \op \sss1{\ee_1}\\
\label{E:Delta2}
&= \sss1\pp \op \flip3\bb' \op 
\sss2{\ee_{\pp+2}-1} \, \sss1{-1} \, \sss2{\ee_{\pp+1}-2} \, \sss1{-1} \, ... \,
\sss1{-1} \, \sss2{\ee_3-2} \, \sss1{-1} \,  \sss2{\ee_2-1} \, \sss1{\ee_1}
\end{align}
Owing to the hypotheses about the exponents~$\ee_\kk$, the
generator~$\sss2{-1}$ does not occur in the expression
of~\eqref{E:Delta2}. If at least one of the inequalities $\ee_{\pp+2} \ge 1$,
$\ee_{\pp+1}, ..., \ee_3 \ge 2$, $\ee_2 \ge1$ is strict, the generator~$\ss2$
occurs in~\eqref{E:Delta2}, and we deduce $\ddd\pp < \bb$. Otherwise,
by definition of normality, $\ee_{\pp+1} = 1$ implies $\bb' = 1$, and
\eqref{E:Delta2} reduces to $\ddd\pp\inv\bb = \sss1{\ee_1}$. If
$\ee_1$ is positive, the main generator~$\ss1$ of~$\sss1{\ee_1}$
occurs positively only, and we have $\ddd\pp < \bb$ again. Finally,
$\ee_1 = 0$ corresponds to $\bb = \ddd\pp$, so $\ddd\pp \le \bb$
holds in all cases.
\end{proof}

We easily deduce the following connection between the braid ordering and
the so-called $\ShortLex$-ordering of the corresponding exponent
sequences.

\begin{prop}
\label{P:Order3}
Assume that $\bb, \bb'$ belong to~$\BB3$. Let $(\ee_\pp, ..., \ee_1)$ and
$(\ee'_\qq, ..., \ee'_1)$ be the exponent sequences of~$\bb$ and~$\bb'$,
respectively. Then $\bb < \bb'$ holds if and only if the sequence
$(\ee_\pp, ..., \ee_1)$ is $\ShortLex$-smaller than the sequence
$(\ee'_\qq, ..., \ee'_1)$, meaning that either $\pp <\qq$ holds, or
we have $\pp = \qq$ and there exists~$\rr$ such that $\ee_\kk =
\ee'_\kk$ holds for $\kk > \rr$ and we have $\ee_\rr < \ee'_\rr$.
\end{prop}

\begin{proof}
Lemma~\ref{L:Order} shows that $\pp< \qq$ implies $\bb < \bb'$. Assume
now $\pp = \qq$, and we have $\ee_\kk = \ee'_\kk$ for $\kk > \rr$ and
$\ee_\rr < \ee'_\rr$. Let $\bb_1 = \sss{\cl{\rr-1}}{\ee_{\rr-1}} ...
\sss2{\ee_2} \sss1{\ee_1}$ and $\bb'_1 = \sss{\cl{\rr}}{\ee'_\rr - \ee_\rr}
\sss{\cl{\rr-1}}{\ee'_{\rr-1}} ...
\sss2{\ee'_2} \sss1{\ee'_1}$. By hypothesis, we have $\bb = \bb_0 \bb_1$
and $\bb' = \bb_0 \bb'_1$ for some braid~$\bb_0$. Truncating
an exponent sequence on the left preserves the normality conditions, hence the
exponent sequence of~$\bb_1$ is $(\ee_{\rr-1}, ..., \ee_1)$, and that
of~$\bb'_1$ is $(\ee'_\rr - \ee_\rr, \ee'_{\rr-1}, ..., \ee'_1)$.
So $\bb_1$ has breadth~$\rr-1$, while $\bb'_1$ has breadth~$\rr$. Then  
Lemma~\ref{L:Order} implies $\bb_1 < \bb'_1$, and $\bb < \bb'$
immediately follows.
\end{proof}

For instance, the definition gives $\ddd\pp = \sss{\cl{\pp+2}}1
\sss{\cl{\pp+1}}2... \sss12 \sss21 \sss10$, implying that the exponent
sequence of~$\ddd\pp$ is $(1, 2, ..., 2, 1, 0)$, with $2$ repeated
$\pp-1$~times: as the latter sequence is $\ShortLex$-minimal among all
length~$\pp+2$ sequences satisfying the normality conditions, we see that
$\ddd\pp$ is indeed the least upper bound of all braids with breadth at
most~$\pp+1$, as stated in Lemma~\ref{L:Order}.

\begin{rema}
The computations above actually reprove that any two braids in~$\BB3$ are 
comparable with respect to the relation~$<$. Indeed, what
Lemma~\ref{L:Order} and Proposition~\ref{P:Order3} prove is that, if the
exponent sequence of~$\bb$ is $\ShortLex$-smaller than that of~$\bb'$,
then the quotient braid $\bb\inv \bb'$ admits at least one expression in
which the generator with highest index occurs only positively.
\end{rema}

\section{Long descending sequences in~$\BB3$}
\label{S:Sequence3}

As the restriction of the braid order to the monoid~$\BB3$ is a well-order, it
admits no infinite descending sequence. However, as the order-type of the
latter well-order is the ordinal~$\om^\om$, there will exist long finite descending
sequences. Here we investigate a powerful method for constructing such
descending sequences by iterating a simple inductive process. The fact that
the sequences are long comes from their apparently growing at each step.
Then, our results are based on the fact that, on the one hand, the well-foundedness
of the braid order forces the sequences to be finite, while, on the other hand, the
sequences are so long that their finiteness cannot be proved in a weak system
like~$\IS1$.

\subsection{$\G_3$-sequences}

The principle is to start with an arbitrary braid in~$\BB3$ and to repeat
some braid transformation until, if ever, the trivial braid (the one with no twist) is
obtained. The transformation at step~$\tt$ consists in removing
one crossing in the considered braid, but then, in all cases but one, reintroducing
$\tt$~new crossings. Thus, the definition is reminiscent of Kirby-Paris' Hydra
Game~\cite{KiP}, with Hercules chopping off one head of the
Hydra and the Hydra sprouting $\tt$~families of new heads. The paradoxical result is
that, contrary to what examples suggest, one always reaches the trivial braid
after finitely many steps.

Our sequences will be defined in terms of the $\flip{}$-normal form of
Definition~\ref{D:Normal3}, and we need some terminology. First, by definition, each block in a
$\flip{}$-normal word, except possibly the leftmost one, has its size at least equal to the minimal
legal size~$\emin_\kk$ introduced in Definition~\ref{D:Normal3}. Our aim will be to remove
crossings in a braid diagram trying to preserve $\flip{}$-normality as much as possible. Therefore,
we are naturally led to considering the blocks whose size strictly exceeds the minimal legal value.

\begin{figure} [htb]
\begin{picture}(112,28)
\put(-0.5,5){\includegraphics{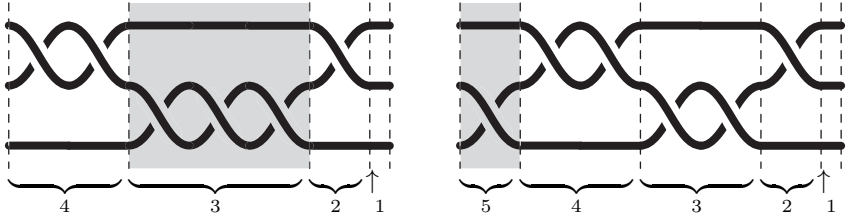}}
\put(0,4){$\underbrace{\hbox to 15mm{}}_4$}
\put(16,4){$\underbrace{\hbox to 23mm{}}_3$}
\put(40,4){$\underbrace{\hbox to 7mm{}}_2$}
\put(47.5,2.5){$\uparrow$}
\put(48.7,-.1){$_1$}
\put(60,4){$\underbrace{\hbox to 7mm{}}_5$}
\put(68,4){$\underbrace{\hbox to 15mm{}}_4$}
\put(84,4){$\underbrace{\hbox to 15mm{}}_3$}
\put(100,4){$\underbrace{\hbox to 7mm{}}_2$}
\put(107.5,2.5){$\uparrow$}
\put(108.7,-.1){$_1$}
\end{picture}
\caption{\smaller\sf Critical position of a braid: in the $\flip{}$-normal 
representative diagram, it corresponds to the rightmost block whose size strictly exceeds the
minimal legal size, if it exists, and to the leftmost block otherwise; on the left,
the example of $\sss22\sss13\ss2$: the exponent sequence is $(2,3,1,0)$, so the critical position
is~$3$, because, with the notation of Definition~\ref{D:Normal3}, we have $\ee_1 = 0 =
\emin_1$, $\ee_2 = 1 = \emin_2$, but $\ee_3 = 3 > \emin_3$; on the right,
the example of $\ss1\sss22\sss12\ss2$: the exponent sequence is $(1,2,2,1,0)$, and the critical
position is the breadth~$5$, because no block has a size exceeding the minimal size.}
\label{F:Critical}
\end{figure}

\begin{defi}
\label{D:Critical3}
(Figure~\ref{F:Critical}) Let $\bb$ be a nontrivial braid in~$\BB3$, and let
$(\ee_\pp, ..., \ee_1)$ be its exponent sequence. The least number~$\rr < \pp$ for which
$\ee_\rr > \emin_\rr$ holds, if such a number exists, or~$\pp$ otherwise, is
called its {\it critical position} in~$\bb$. 
\end{defi}

Thus, the critical position of~$\bb$ corresponds to the rightmost block in the $\flip{}$-normal
expression of~$\bb$ whose size is not minimal---hence the rightmost block in which one
can remove one crossing without destroying normality---if such a block exists,
and to the leftmost block otherwise. For instance, the critical position in~$\D_3$ is~$1$, as the
length of the final block of~$\ss1$'s, here~$1$, is positive, hence strictly larger than the
minimal value~$\emin_1$. 

We are ready to define $\G_3$-sequences. In order to apply the principles described above, in
particular to preserve $\flip{}$-normality, we have to choose a position where to remove one
crossing, and this is where we use the critical position.

\begin{figure} [htb]
\begin{picture}(108,27)(0,0)
\put(0,5){\includegraphics{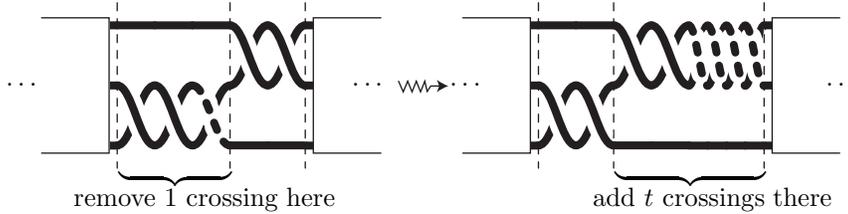}}
\put(12,5){$\underbrace{\hbox to15mm{}}$}
\put(6,0){remove $1$ crossing here}
\put(78,5){$\underbrace{\hbox to20mm{}}$}
\put(75,0){add $\tt$ crossings there}
\put(-3,16){\dots}
\put(43,16){\dots}
\put(56,16){\dots}
\put(106,16){\dots}
\end{picture}
\caption{\smaller\sf Inductive construction of the $\G_3$-sequence: at step~$\tt$---here $\tt=4$---we remove one crossing in the critical block, but $\tt$~new crossings appear in the next block, if it exists,  \ie, if the critical block is not the final block of~$\ss1$'s.}
\label{F:Rule}
\end{figure}

\begin{defi}
\label{D:Game}
(Figure~\ref{F:Rule}) Assume that $\bb$ is a nontrivial $3$-braid, and $\tt$ is a natural number. Let $\ww$ be the $\flip{}$-normal word representing~$\bb$, and~$\rr$ be the critical position in~$\bb$. Then we define $\ww\OP\tt$ to be the word obtained from~$\ww$ by removing one letter in the $\rr$th block, and adding $\tt$~letters in the $(\rr-1)$th block if the latter exists, \ie, if $\rr \ge 2$ holds. We define $\bb\OP\tt$ to be the braid represented by~$\ww\OP\tt$, and the {\it $\G_3$-sequence from~$\bb$} to be the sequence $(\bb_0, \bb_1, ...)$ defined by $\bb_0 = \bb$ and $\bb_\tt = \bb_{\tt-1}\OP\tt$; the sequence stops when the trivial braid~$1$ is possibly obtained.
\end{defi}

\begin{exam}
\label{X:Battle}
Let $\bb = \sss22\sss12$. The $\G_3$-sequence from~$\bb$ is as follows:
$$(\sss22\sss12, 
\sss22\ss1,
\sss22,
\ss2\sss13,
\ss2\sss12,
\ss2\ss1,
\ss2,
\sss17, 
\sss16, 
\sss15, 
\sss14, 
\sss13, 
\sss12, 
\ss1, 
1),$$
\ie, in this case, we reach the trivial braid in 14~steps.
Not all examples are so simple. The reader can check that,
starting from~$\D_3$, \ie, $\ss1 \ss2 \ss1$, one reaches the trivial braid after
$30$~steps, whereas, starting from $\sss12\sss22\sss12$, a braid with six
crossings only, one does reach the trivial braid after no less than
$90,159,953,477,630$~steps...
\end{exam}

\begin{rema}
In principle, constructing a $\G_3$-sequence entails finding at each step the
$\flip{}$-normal word that represents the current braid. Actually, this procedure
has to be performed at the initial step only. Indeed, our definition of the critical
position guarantees that, if $\ww$ is a $\flip{}$-normal word, then, for
every~$\tt$, the word~$\ww\OP\tt$ is $\flip{}$-normal as well. So, if $\ww$ is
the $\flip{}$-normal word representing~$\bb$, then, for each~$\tt$, the
$\flip{}$-normal word representing $\bb\OP1\OP2...\OP\tt$ is
$\ww\OP1\OP2...\OP\tt$. In other words, provided we start with a
$\flip{}$-normal diagram, we can play with braid diagrams without worrying
about normalization.
\end{rema}

\subsection{Finiteness of $\G_3$-sequences}

The first result is that, although very long $\G_3$-sequences exist, no
such sequence is infinite, \ie, Proposition~A of the introduction.

\begin{prop}
\label{P:Termination}
For each $3$-braid~$\bb$, the $\G_3$-sequence from~$\bb$ is finite,
\ie, there exists a finite number~$\tt$ satisfying $\bb\OP1\OP2...\OP\tt =
1$.
\end{prop}

Proposition~\ref{P:Termination} directly follows from the conjunction of
two results, namely that, according to Theorem~\ref{T:Order}$(ii)$, the
braid order on~$\BB3$ is a well-order, hence possesses no infinite
descending sequence, and that every $\G_3$-sequence is descending with
respect to that order. The latter result is a direct consequence of

\begin{lemm}
\label{L:Descending3}
For each braid~$\bb$ in~$\BB3$ and each number~$\tt$, we have $\bb >
\bb\OP\tt$.
\end{lemm}

\begin{proof}
Assume $\bb' = \bb\OP\tt$. There are three possible cases, according to the 
critical position~$\rr$ in~$\bb$. For $\rr = 1$, \ie, if the critical block is the
final block of~$\ss1$'s, then we directly have ${\bb'}\inv\bb = \ss1$, hence
$\bb' < \bb$. For $\pp-1 \ge \rr \ge2$, and for $\pp = \rr \ge 2$ with
$\ee_\pp \ge 2$, the exponent sequence of~$\bb'$ is obtained from that
of~$\bb$ by replacing some subsequence $(\ee_\rr, \ee_{\rr-1})$ with
$(\ee_\rr-1, \ee_{\rr-1} + \tt)$. Hence the exponent sequence of~$\bb'$ is
$\ShortLex$-smaller than that of~$\bb$, and Proposition~\ref{P:Order3}
implies $\bb' < \bb$. Finally, for $\pp = \rr \ge 2$ with $\ee_\pp =1$, the
exponent sequence of~$\bb'$ is obtained from that of~$\bb$ by replacing
some subsequence $(1, \ee_{\rr-1})$ with $(\ee_{\rr-1} + \tt)$, and, again,
the exponent sequence of~$\bb'$ is $\ShortLex$-smaller than that of~$\bb$,
and Proposition~\ref{P:Order3} implies $\bb' < \bb$.
\end{proof}

\begin{rema}
Some variants are possible in the definition of $\G_3$-sequences.  In
particular, $\G_3$-sequences are deterministic: for each non-trivial
braid~$\bb$ and each number~$\tt$, the braid~$\bb\OP\tt$ is uniquely defined. Actually, we
could instead consider at each step an arbitrary permitted block rather than the critical block, a
block being called permitted whenever its size exceeds the minimal legal size of
Definition~\ref{D:Normal3}---so that the critical block is just the rightmost permitted block. In this
way, we obtain in general many sequences from an initial braid~$\bb$. However, the argument of
Lemma~\ref{L:Descending3} remains valid, so each such sequence has to be finite. Such a variant
can be described as a game (or a battle) against a braid~$\bb$: the player tries to destroy the
braid, namely to reduce it to the trivial braid; the rule is that, at step~$\tt$, the player chooses
one permitted position and removes one crossing from the corresponding block; then, unless the
block was the final block of~$\ss1$'s, the (nasty!) braid lets $\tt$~new crossings appear in
the next block. Lemma~\ref{L:Descending3} guarantees that every
battle against every $3$-braid is won---as is every battle against every hydra
in~\cite{KiP, Car}.
\end{rema}

\subsection{An unprovability result}

We turn to Theorem~A of the introduction, \ie, we prove
that the finiteness of $\G_3$-sequences cannot be proved from the axioms
of~$\IS1$---see Appendix for a brief introduction to that system. 

\begin{thrm}
\label{T:UnprovGame}
Proposition~\ref{P:Termination} is an arithmetic
statement that cannot be proved in the system~$\IS1$.
\end{thrm}

Theorem~\ref{T:UnprovGame} follows from the fact that the length of 
well-chosen $\G_3$-sequences grows faster than any recursive function whose
totality is provable from~$\IS1$ axioms\footnote{\ie, every function~$\ff$ such that
$\yy = \ff(\xx)$ can be expressed by some $\Sigma_1$-formula~$\For(\xx,\yy)$ such that
$\forall\xx \exists\yy (\For(\xx,\yy))$ is provable from~$\IS1$ axioms.}. 

\begin{defi}
For each braid~$\bb$, we denote by~$\TT(\bb)$ the length  of the
$\G_3$-sequence from~$\bb$, \ie, the smallest integer~$\tt$ satisfying
$\bb\OP1\OP2...\OP\tt = 1$. 
\end{defi}

It is well-known that the Ackermann function eventually dominates, provably in~$\IS1$,  all
functions that are  provably total in $\IS1$. Therefore, in order to prove
Theorem~\ref{T:UnprovGame}, it is enough to find an explicit\footnote{more
precisely, a primitive recursive one since the argument has to take place inside~$\IS1$} sequence
of braids~$(\bb_0, \bb_1, ...)$ such that the function $\xx \mapsto \TT(\bb_\xx)$ eventually
dominates the Ackermann function. In order to do that, we shall resort to the so-called
fundamental sequences of ordinals and the Hardy hierarchy of fast growing functions~\cite{BCW}.

Fundamental sequences of ordinals are obtained by selecting, for each (limit) ordinal~$\l$
in an interval, here $[0, \om^\om]$, a distinguished increasing sequence cofinal in~$\l$. To this aim,
referring to the Cantor Normal Form of ordinals---see Appendix---we write $\a \eqCNF \b + \om^\delta$ to
mean that $\a = \b + \om^\delta$ holds and, in addition, $\om^\delta$ is the last factor in the Cantor
Normal Form of~$\a$, \ie, $\b$ is either~$0$, or it can be written as $\om^{\b_1} + ... + \om^{\b_\mm}$
with $\b_\mm \ge \delta$. Then we observe that every limit ordinal~$\l$ below~$\om^\om$ can be
uniquely expressed as $\l \eqCNF \g + \om^{\rr+1}$. 

\begin{defi}
\label{D:Fundamental3}
For $\l$ a limit ordinal below~$\om^\om$, say $\l \eqCNF \g + \om^{\rr+1}$, and $\xx$ is a nonnegative
integer, we define $\l[\xx] = \g + \om^\rr \cdot
\xx$. Moreover, we put $\om^\om[\xx] = \om^\xx$. 
\end{defi}

For technical convenience, the definition is extended to non-limit ordinals by setting $0[\xx]= 0$
and $(\alpha+1)[\xx]=\alpha$ for every~$\xx$.

By construction, for each limit ordinal~$\l$, the sequence $\l[0], \l[1], ...$ is
increasing and cofinal in~$\l$.  The main technical result we use consists
in associating with every $3$-braid an ordinal below~$\om^\om$ so that there exists a simple
connection between the operations $\bb \mapsto \bb\OP\tt$ in the
braid monoid~$\BB3$ and $\b \mapsto \b[\tt]$ in the ordinal interval~$[0,
\om^\om)$.

\begin{defi}
For $\bb$ a $3$-braid with exponent sequence $(\ee_\pp, ..., \ee_1)$, we put
\begin{equation}
\label{E:Ordinal3}
\ord(\bb) = \om^{\pp-1} \cdot \ee_\pp + \sum_{\pp>\kk\ge1} \om^{\kk-1} \cdot 
(\ee_\kk - \emin_\kk).
\end{equation} 
\end{defi}

The idea of the definition is simply to measure by which amount the exponent
sequence of~$\bb$ exceeds the minimal legal values. 

\begin{exam}
\label{X:Ord}
The exponent sequence of the braid~$\ddd\pp$ introduced in
Lemma~\ref{L:Order} is $(1, 2, ..., 2, 1, 0)$, with $2$ repeated~$\pp-1$
times. We deduce $\ord(\ddd\pp) = \om^{\pp+1}$ for~$\pp \ge 1$.
Similarly, we noted in the proof of Lemma~\ref{L:Order} that
$\D_3^\pp = \ddd\pp\sss1\pp$ holds, and we deduce $\ord(\D_3^\pp) =
\om^{\pp+1} + \pp$ for~$\pp \ge 1$.
\end{exam}

It can be checked that, for each braid~$\bb$ in~$\BB3$, the
ordinal~$\ord(\bb)$ is the rank of~$\bb$ in the well-ordering~$(\BB3,<)$,
but we shall not use this result. Note that the expression of~$\ord(\bb)$
given in~\eqref{E:Ordinal3} is always in Cantor Normal Form---see
Appendix. 

\begin{lemm}
\label{L:OneConnection} 
For every nontrivial $\bb$ in~$\BB3$ and every $\tt$ in~$\Nat$, we have
\begin{equation}
\label{E:Connection}
\ord(\bb\OP\tt) = \ord(\bb)[\tt].
\end{equation}
\end{lemm}

\begin{proof}
Put $\bb' = \bb\OP\tt$, and let $\ww = \sss{\cl\pp}{e_\pp} ... \sss2{\ee_2} \sss1{\ee_1}$ be the $\flip{}$-normal 
expression of~$\bb$. We consider the various possible values of the
critical position~$\rr$ of~$\bb$. .

Case 1: $\rr = 1$. In this case, we have $\ee_1 \ge 1$, whence
$\bb' = \sss{\cl\pp}{e_\pp} \, ... \,  \sss2{\ee_2} \sss1{\ee_1-1}$.
By definition of the function~$\ord$, we obtain, both for $\pp = 1$ and $\pp \ge2$,
$$\ord(\bb) = \ord(\bb') + 1.$$
As~$\ord(\bb)$ is a successor ordinal, the latter equality also reads $\ord(\bb') = \ord(\bb)[\tt]$.

Case 2: $\pp > \rr \ge 2$. Then, we have $\bb' = \sss{\cl\pp}{e_\pp} \, ... \,
\sss{\cl{\rr+1}}{\ee_{\rr+1}} 
\sss{\cl\rr}{\ee_\rr-1} 
\sss{\cl{\rr-1}}{\ee_{\rr-1}+\tt} 
\sss{\cl{\rr-2}}{\ee_{\rr-2}} 
\,...\,
 \sss1{\ee_1}$, so, taking into account the hypothesis that $\kk < \rr$ implies $\ee_\kk = \emin_\kk$, we deduce
\begin{gather*}
\ord(\bb) = \om^{\pp-1} \cdot \ee_\pp + \sum_{\pp > \kk > \rr} \om^{\kk-1}\cdot 
(\ee_\kk - \emin_\kk) + \om^{\rr-1}\cdot 
(\ee_\rr - \emin_\rr),\\
\ord(\bb') = \om^{\pp-1} \cdot \ee_\pp + \sum_{\pp > \kk > \rr} \om^{\kk-1}\cdot 
(\ee_\kk - \emin_\kk) + \om^{\rr-1}\cdot 
(\ee_\rr - \emin_\rr-1) + \om^{\rr-2}\cdot \tt. 
\end{gather*}
Putting $\g = \om^{\pp-1} \cdot \ee_\pp + \sum_{\pp > \kk > \rr} \om^{\kk-1}\cdot  (\ee_\kk - \emin_\kk) + \om^{\rr-1}\cdot  (\ee_\rr - \emin_\rr-1)$, the latter values read $\ord(\bb) = \g + \om^{\rr-1}$ and $\ord(\bb') = \g + \om^{\rr-2} \cdot \tt$, so $\ord(\bb') = \ord(\bb)[\tt]$ holds.

Case 3: $\pp = \rr \ge 2$. In this case, we find $\bb' = \sss{\cl\pp}{\ee_\pp-1} 
\sss{\cl{\pp-1}}{\ee_{\pp-1}+\tt} 
\sss{\cl{\pp-2}}{\ee_{\pp-2}} 
\,...\,
 \sss1{\ee_1}$. As in Case~$2$, we have $\ee_{\pp-1} = \emin_{\pp-1}$, and we obtain
\begin{gather*}
\ord(\bb) = \om^{\pp-1} \cdot \ee_\pp,
\text{\quad and \quad}
\ord(\bb') = \om^{\pp-1} \cdot (\ee_\pp-1) + \om^{\pp-2}\cdot \tt, 
\end{gather*}
\ie, again, $\ord(\bb) = \g + \om^{\rr-1}$ and $\ord(\bb') =  \g +
\om^{\rr-2} \cdot \tt$ when we put $\g = \om^{\pp-1} \cdot (\ee_\pp-1)$.
So $\ord(\bb') = \ord(\bb)[\tt]$ holds in this case as well.
\end{proof}

We easily deduce a comparison between the function~$\TT$ measuring the length of $\G_3$-sequences and the functions~$H_\a$ of the Hardy hierarchy. We recall the definition of the latter.

\begin{defi}
\label{D:Hardy3}
For $\a \le \om^\om$, the functions $H_\a : \Nat \to \Nat$ are defined by 
\begin{equation}
\label{E:Hardy3}
H_{\alpha}(\xx):=
\begin{cases}
\xx & \mbox{ if } \alpha = 0,\\
H_{\beta}(\xx+1) & \mbox{ if } \alpha = \beta+1,\\
H_{\alpha[\xx]}(\xx+1) & \mbox{ if } \alpha \mbox{ is a limit ordinal}.
\end{cases}
\end{equation}
\end{defi}

For instance, we have $H_\rr(\xx) = \xx+\rr$ for each natural number~$\rr$, then $H_\om(\xx) = 2\xx+1$, $H_{\om+\rr}(\xx) = 2\xx + 2\rr + 1$, $H_{\om\cdot2}(\xx) = 4\xx + 3$, etc.
It is known---see for instance~\cite{BCW}---that the function~$H_{\om^\om}$ is ackermannian, \ie,
it is a slight variant of  the Ackermann function. 

An easy induction from the definition---see~\cite{BCW} again---gives for every~$\b \le \om^\om$
and for every~$\kk$
\begin{equation}
\label{E:Hardy}
H_\b(\kk) = \min\{\tt \mid \b[\kk]...[\kk+\tt-1]=0\} + \kk.
\end{equation}
Then we obtain the main comparison result:

\begin{prop}
\label{P:Hardy} 
Let $\bb$ be a $3$-braid with $\ord(\bb) = \b$. Then, for each~$\kk$, we have
\begin{equation}
\label{E:HardyConnection}
\TT(\bb\sss1{\kk}) = H_\b(\kk+1) - 1.
\end{equation}
\end{prop}

\begin{proof}
By construction, we have $\bb\sss1{\kk}\OP1\OP2...\OP{\kk} = \bb$.
So $\bb\OP{\kk+1} ...\OP{\kk+\tt} = 1$ is equivalent to
$\bb\sss1{\kk} \OP1...\OP{\kk+\tt} = 1$, and, therefore,
$\TT(\bb\sss1{\kk})$ equals $\kk$ plus the smallest~$\tt$ for which
$\bb\OP{\kk+1} ... \OP{\kk+\tt} = 1$ holds. For all $\kk$ and~$\tt$,
repeated applications of Lemma~\ref{L:OneConnection} yield:
\begin{equation*}
\ord(\bb\OP{\kk+1} ...  \OP{\kk+\tt})
= \ord(\bb)[\kk+1] ... [\kk+\tt].
\end{equation*}
Now $\bb\OP{\kk+1} ...  \OP{\kk+\tt}$ is the trivial braid~$1$ if and
only if the associated ordinal is~$0$, and, therefore, the smallest~$\tt$ for
which
$\bb\OP{\kk+1} ...  \OP{\kk+\tt} = 1$ holds is the smallest~$\tt$ for
which
$\ord(\bb)[\kk+1] ... [\kk+\tt] = 0$ holds, so we obtain
\begin{equation*}
\TT(\bb\sss1{\kk}) 
= \kk + \min\{\tt \mid \b[\kk+1]...[\kk+\tt]=0\},
\end{equation*}
which, by~\eqref{E:Hardy}, is~$H_\b(\kk+1) - 1$.
\end{proof}

We are now ready to complete the main argument.

\begin{proof}[Proof of Theorem~\ref{T:UnprovGame}]
Define~$\UU$ by $\UU(0) = 2$, $\UU(1) = 5$, and $\UU(\kk) =  \TT(\D_3^{\kk-1} \ss1) +1$ for
$\kk \ge 2$. By~\eqref{E:Deltadelta}, we have $\D_3^{\kk} = \ddd{\kk-1}
\sss1{\kk}$ for $\kk \ge 1$, with $\ord(\ddd0) = 0$ and $\ord(\ddd{\kk-1}) = \om^\kk$ for $\kk \ge 2$.
So, for $\kk \ge 2$, \eqref{E:HardyConnection} plus the definition of the function~$H_{\om^\om}$ give 
$$\UU(\kk) = \TT(\ddd{\kk-1} \sss1{\kk}) + 1 = 
H_{\om^\kk}(\kk+1) = H_{\om^\om}(\kk)$$
---actually $\UU(\kk) = H_{\om^\om}(\kk)$ holds for each~$\kk$ owing to the values at~$0$ and~$1$.  
So, the function~$\UU$ is the ackermannian function~$H_{\om^\om}$, and, therefore, it cannot
be primitive recursive. Now, if the finiteness of
$\G_3$-sequences were provable in~$\IS1$, the  function~$\UU$ would be
provably total in $\IS1$---see Appendix.  But every provably total function
of~$\IS1$ is primitive recursive~\cite{Min, Par}---see for instance~\cite{FaW}.
\end{proof}

\section{Friedman-style results and phase transitions}
\label{S:Transition}

With $\G_3$-sequences, we considered descending sequences of a particular type. We shall 
now consider more general sequences, where the entries no longer obey a
particular formation law, but only satisfy some growth conditions defined in
terms of Garside's complexity, in the spirit of the 
sentences considered by H.\,Friedman in~\cite{Fri}. The main result here is that
there exists a precise description of the conditions that lead from
$\IS1$-provability to $\IS1$-unprovability, thus witnessing for a quick phase
transition phenomenon analogous to those investigated in~\cite{Wei2, Wei3,
Wei4}.

\subsection{The complexity of a $3$-braid}
\label{S:Complexity}

In the sequel, we need some measure for the complexity of a braid. We shall
resort to the most usual such measure, namely the complexity, or canonical
length, derived from Garside's theory---see for instance \cite[Chapter~9]{Eps}. 

\begin{defi}
\label{D:Canonical}
We say that a positive $3$-strand braid~$\bb$ has {\it complexity}~$\ell$, written $\can\bb = \ell$, if $\bb$ is a left divisor
of~$\D_3^\ell$, \ie, there exists a braid~$\bb'$ in~$\BB3$ satisfying
$\bb\bb' = \D_3^\ell$, and $\ell$ is minimal with that property.
\end{defi}

If $\bb$ is a braid, we use $\lg\bb$ for the common length of all braid words
representing~$\bb$. By Garside theory, $\lg\bb \le \ell$ implies $\can\bb
\le \ell$; on the other hand, $\bb$ dividing~$\D_3^\ell$ implies 
$\lg\bb \le \lg{\D_3^\ell} = 3\ell$, so, for each positive $3$-strand
braid~$\bb$, we have the inequalities
\begin{equation}
\label{E:Comparison}
\can\bb \le \lg\bb \le 3\can\bb,
\end{equation}
which show that, in the case of~$\BB3$, the metrics associated with the length
and the complexity are quasi-isometric. Although all arguments below could
be completed using the inequalities of~\eqref{E:Comparison} exclusively---at
the expense of modifying some parameters---it will be convenient to resort to
a more precise connection between the length, the complexity, and
the breadth of a $3$-braid. This connection relies on a simple relation
between the greedy normal form and the $\flip{}$-normal form of a
$3$-braid that was first observed by J.\,Mairesse, and that is of independent
interest.

\begin{lemm}
\label{L:CanLength}
$(i)$ For $\bb$ a $3$-braid, let $\expo\bb$ denote the maximal
integer~$\dd$ such that $\D_3^\dd$ is a divisor of~$\bb$. Then $\expo\bb$
is the maximal~$\dd$ such that the exponent sequence of~$\bb$ has
the form $(..., \ee_\dd, 2, ..., 2, 1, \ee_1)$ with $\ee_\dd \ge 1$ and
$\ee_1 \ge \dd$.

$(ii)$ For each $3$-braid~$\bb$, we have
\begin{equation}
\label{E:ComparisonBis}
\can\bb = \lg\bb - \pp - \expo\bb + C,
\end{equation}
where $\pp$ is the breadth of~$\bb$ and $C$ is $0$, $1$, or~$2$.
\end{lemm} 

\begin{proof}
It is standard---see for instance \cite[Chap.~9]{Eps} or the
introduction of~\cite{Dhh}---that every
$3$-braid admits a unique expression of the form $\ww_\rr ... \ww_1
(\ss1\ss2\ss1)^\dd$, where, for each~$\kk$, the word $\ww_\kk$ is either
$\ss1$, or $\ss2$, or $\ss1\ss2$, or $\ss2\ss1$, and the last letter
of~$\ww_{\kk+1}$ is the first letter of~$\ww_\kk$. This expression is called
the right greedy (or Garside) normal form of~$\bb$, and we have then
$\can\bb = \rr + \dd$. By grouping the letters, the greedy normal form
of~$\bb$ can be uniquely written as $\ww = \sss{\cl\qq}{\dd_\qq} ...
\sss2{\dd_2}
\sss1{\dd_1} (\ss1\ss2\ss1)^\dd$ with $\qq \ge 0$, $\dd_\qq, ..., \dd_2 \ge
1$ and $\dd_1 \ge 0$, and the above formula gives
\begin{equation}
\label{E:CanSeq}
\can\bb = 
\begin{cases}
\dd&\mbox{for $\qq = 0$},\\
\dd_\qq + \cdots + \dd_1 + \dd - \qq + 1
&\mbox{for $\qq > 0$ with $\dd_1 > 0$},\\
\dd_\qq + \cdots + \dd_1 + \dd - \qq + 2
&\mbox{for $\qq > 0$ with $\dd_1 = 0$}.
\end{cases}
\end{equation}
An easy computation shows that the exponent sequence
of~$\D_3^\dd$ is $(1,2^{(\dd-1)}, 1, \dd)$, where $2^{(\mm)}$
stands for $2, ..., 2$ with~$2$ repeated $\mm$~times. As conjugating by~$\D_3$, \ie, applying the flip
automorphism~$\flip3$, exchanges~$\ss1$ and~$\ss2$, we have 
$\sss1{\dd_1}\D_3^\dd = \D_3^\dd\sss1{\dd_1}$ if $\dd$ is
even, and
$\sss2{\dd_2}\D_3^\dd = \D_3^\dd\sss1{\dd_2}$ is $\dd$ if
odd, we obtain that $\bb$ is represented by the word whose
exponent sequence is
\begin{equation}
\label{E:ExpSeq}
\begin{cases}
(\dd_\qq, ..., \dd_3, \dd_2+1, 2^{(\dd-1)}, 1, \dd+\dd_1)
&\mbox{for $\dd$ even},\\
(\dd_\qq, ..., \dd_3, \dd_2, \dd_1+1, 2^{(\dd-1)}, 1, \dd)
&\mbox{for $\dd$ odd with $\dd_1 > 0$},\\
(\dd_\qq, ..., \dd_4, \dd_3+1, 2^{(\dd-1)}, 1, \dd+\dd_2)
&\mbox{for $\dd$ odd with $\dd_1 = 0$}.
\end{cases}
\end{equation}
In each case, the above sequence satisfies the requirements of
Proposition~\ref{P:Normal3}, hence, by uniqueness, it is the exponent
sequence of~$\bb$.

Now, the explicit form of the sequences occurring in~\eqref{E:ExpSeq} shows
that, in each case, the parameter~$\dd$ corresponds to the longest suffix of
the form $(\ee_\dd, 2^{(\dd-1)},1,\ee_1)$ with $\ee_\dd \ge 1$ and $\ee_1
\ge \dd$, which gives~$(i)$. 

Similarly, one easily deduces from~\eqref{E:ExpSeq} that the breadth
of~$\bb$ is $\qq + \dd + C'$, with $C' = 0$ with $\dd$ is even, and $C' = 1$
(\resp $-1$) when $\dd$ is odd with $\dd_1>0$ (\resp $= 0$).
Plugging these values in~\eqref{E:CanSeq} and using the relation 
$\lg\bb = \dd_\qq + \cdots + \dd_1 + 3 \dd$ gives~\eqref{E:ComparisonBis}
with $C = 0$ for $\qq = 0$, $C=1$ for $\dd$ even with $\dd_1>0$ and
$\dd$ odd with $\dd_1 = 0$, and $C=2$ for $\dd$ even with $\dd_1=0$ and
$\dd$ odd with $\dd_1 > 0$.
\end{proof}

\subsection{Combinatorial well-foundedness of the braid ordering}

The idea is to consider combinatorial principles that capture some finitary
aspects of the well-foundedness of the braid ordering. The first obvious
observation is that, for each constant~$\kk$, there exist only finitely many
braids with complexity bounded by~$\kk$, and, therefore, there exists
an obvious upper bound on the length of possible decreasing sequences
consisting of such braids.

\begin{prop}
\label{P:Constant}
For each~$\kk$, there exists~$\mm$ such that there exists no descending
sequence $(\bb_0, ..., \bb_\mm)$ in~$\BB3$ such that $\can{\bb_\ttii} \le \kk$
holds for each~$\ttii$.
\end{prop}

\begin{proof}
By~\eqref{E:Comparison}, $\can\bb \le \kk$ implies $\lg\bb \le 3\kk$,
so there are at most $1 + 2 + 4 + {\cdot}{\cdot}{\cdot} + 2^{3\kk}$~braids of complexity bounded
by~$\kk$. Hence, by the pigeonhole principle, the expected result holds with $\mm = 2^{3\kk+1}$.
\end{proof}

We now relax the condition that the complexity of the braids is
bounded by a fixed number into a weaker condition involving a function parameter.

\begin{defi}
For $\ff: \Nat \to \Nat$, a sequence of braids $(\bb_0, ..., \bb_\mm)$ is said to
be {\it $(\kk,\!\ff)$-simple} if, for each~$\ttii$, we have $\can{\bb_\ttii} \le
\kk + \ff(\ttii)$. We denote by~$\WO\ff$ the principle (``Well-Order
Property of
$(\BB3, <)$ w{.}r{.}t{.}~$\ff$''):
\begin{quote}
For each~$\kk$, there exists~$\mm$ such that there is no
$(\kk,\!\ff)$-simple descending sequence of length~$\mm$ in~$(\BB3,<)$.
\end{quote}
\end{defi}

Roughly speaking, $\WO\ff$ says that there is no very long
descending sequence of braids with a complexity bounded by~$\ff$.
Formally, it is expressed as
$$
\forall\kk \, \exists\mm \, \forall\bb_0, ..., \bb_\mm \in\BB3\, 
\Big(\forall\ttii\le\mm(\can{\bb_\ttii}\le \kk + \ff(\ttii)) \Rightarrow
\exists\ttii<\mm(\bb_\ttii \not> \bb_{\ttii+1})\Big)\footnote{It can be seen that $\WO\ff$ is---or rather can
be encoded in---a $\Pi_2^0$-statement in the language of arithmetic enriched by a name for~$\ff$, \ie, a
formula of the form $\forall\xx_1 \exists\xx_2 (\For)$ with~$\For$ containing bounded quantifiers only.}.
$$
With this terminology, Proposition~\ref{P:Constant} says that, if $\ff$ is a
constant function, then the principle $\WO\ff$ is true. Actually,  $(\BB3, <)$ being well-ordered
easily---yet non-constructively---implies

\begin{prop}
\label{P:Any}
For each function~$\ff$, the sentence~$\WO\ff$ is true.
\end{prop}

\begin{proof}
We use a compactness argument. For each~$\kk$, let $\TT_\kk$ be the set
of all finite $(\kk,\!\ff)$-simple descending sequences in~$\BB3$. Say that 
$(\bb_1, ..., \bb_\mm) \prec (\bb'_1, ..., \bb'_{\mm'})$ holds if we have
$\mm <\mm'$ and $\bb'_\ttii = \bb_\ttii$ for $\ttii \le \mm$, \ie, if the
latter sequence extends the former. Then $(\TT_\kk, \prec)$ is a
partially ordered set, more precisely a tree, as the predecessors of a
length~$\mm$ sequence consist of its prefixes, and therefore are linearly
ordered by~$\prec$. Now we observe that the tree $(\TT_\kk, \prec)$ is
finitely branching, \ie, a given sequence admits only finitely many
immediate $\prec$-successors. Indeed, by def\-inition of
$(\kk,\!\ff)$-simplicity, the possible successors of a sequence $(\bb_1, ...,
\bb_\mm)$ are of the form $(\bb_1, ..., \bb_\mm, \bb)$ with $\bb$
subject to the constraint $\can\bb \le \kk + \ff(\mm+1)$, and there are
finitely many such braids~$\bb$. On the other hand, the fact that $(\BB3,
{<})$ is a well-ordered set implies that $(\TT_\kk, \prec)$ has no infinite
branch. By K\"onig's Lemma, this implies that $\TT_\kk$ is finite. Hence
there exists a number~$\mm$ such that all sequences in~$\TT_\kk$
have length strictly smaller than~$\mm$. This number~$\mm$ witnesses for
the principle~$\WO\ff$.
\end{proof}

We shall now investigate the logical strength of the principle~$\WO\ff$
when the parameter function~$\ff$ varies. The above easy proof shows that
$\WO\ff$ is true for each~$\ff$, but, as it involves K\"onig's Lemma, it is
not formalizable in a weak system like~$\IS1$. It is somehow
surprising that, for certain natural choices of~$\ff$, the statement is actually
unprovable in~$\IS1$. The most striking results will be established in
Section~\ref{SS:Transition} below. For the moment, we shall establish
that the jump from $\IS1$-provability to $\IS1$-unprovability occurs somewhere between
the constant function and the square function. To this end, we shall use
the following result that controls the complexity of~$\bb \OP \tt$ in
terms of that of~$\bb$.

\begin{lemm}
\label{L:Complexity}
For every $3$-braid~$\bb$ and every number~$\tt$, we have
\begin{equation}
\label{E:Complexity}
\can{\bb \OP \tt} \le \can\bb + \tt + 3.
\end{equation}
\end{lemm}

\begin{proof}
We use the evaluation of the complexity given in
Lemma~\ref{L:CanLength}. First, when we go from~$\bb$ to~$\bb
\OP \tt$, the word length increases by at most~$\tt-1$. Next, by
construction, the breadth of~$\bb\OP\tt$ is either that of~$\bb$, or is that
of~$\bb$ diminished by~$1$. Finally, we have $\expo{\bb \OP \tt} \ge
\expo\bb - 1$. Indeed, by Lemma~\ref{L:CanLength}$(i)$, the
value of~$\dd$ corresponds to the longest suffix of the flip normal form that
has the form $(\ee_\dd, 2^{(\dd-1)}, 1, \ee_1)$ with $\ee_1 \ge \dd$. When
we go from~$\bb$ to~$\bb \OP \tt$, the only case when this longest suffix
can be changed corresponds to the case when $\ee_1$ equals~$\dd$
in~$\bb$, and it becomes $\ee_1-1$ in~$\bb\OP\tt$. In all cases when the
critical block of~$\bb$ is not the rightmost block, the parameter~$\dd$ is
simply~$0$, and it cannot decrease. Taking into account the fact that the
constants~$C$ associated with~$\bb$ and with~$\bb\OP\tt$ can differ by at
most~$2$, we deduce $\can{\bb\OP\tt} \le
\can\bb + (\tt-1) + 1 + 1 + 2$ from~\eqref{E:ComparisonBis}.
\end{proof}

\begin{thrm}
\label{T:Jump}
Let $\Const\rr$ denote the constant function with value~$\rr$, and $\square$
be defined by $\square(\xx) = \xx^2$.

$(i)$ For each~$\rr$, the principle $\WO{\Const\rr}$ is provable
from~$\IS1$.

$(ii)$ The principle $\WO\square$ is not provable
from~$\IS1$.
\end{thrm}

\begin{proof}
$(i)$ The counting argument of Proposition~\ref{P:Constant} goes through
in~$\IS1$---as well as in the weaker system $\IDEXP$.

$(ii)$ Let $\bb$ be an arbitrary positive $3$-strand braid. We prove that the $\G_3$-sequence starting
from~$\bb$ is $(\can\bb + 6, \square)$-simple. Note that the argument below can be done in~$\IS1$.
Let $\kk = \can\bb + 6$. Put $\bb_0 = \bb$, and let $\bb_\ttii$ the $\tt$th entry in
the $\G_3$-sequence from~$\bb$. By Lemma~\ref{L:Complexity}, we obtain
$$\can{\bb_\ttii} \le \can\bb + (1+3) + \cdots + (\ttii + 3) = \can\bb +
\frac12\ttii^2 + \frac72\ttii.$$
For $\ttii$ a nonnegative integer, the latter value is bounded above by
$\ttii^2+6$, so, in each case, we have $\can{\bb_\ttii} \le \can\bb + 6 +
\ttii^2$, \ie,
$\can{\bb_\ttii} \le \kk + \square(\ttii)$. So $(\bb_0, ..., \bb_\mm)$  is $(\kk,\!\square)$-simple. 

Now, assume that the principle $\WO\square$ is provable from~$\IS1$.
Then, for the chosen $k$, one can prove from~$\IS1$ the existence of a
constant~$\mm$ such that every descending $(\kk,\!\square)$-simple
sequence has length less than~$\mm$. So, in particular, the $\G_3$-sequence
$(\bb_0, \bb_1, ...)$   from~$\bb$ cannot be descending for more
than~$\mm$ steps, which means that its length is at most~$\mm$. This being
expressible in~$\IS1$, we conclude that the finiteness of $\G_3$-sequences is
provable from~$\IS1$, contradicting Theorem~\ref{T:UnprovGame}.
\end{proof}

So, at this point, we know that the transition between $\IS1$-provability and
$\IS1$-unprovability for~$\WO\ff$ occurs somewhere between
constant functions and the square function---as illustrated in
Figure~\ref{F:Threshold}. We shall improve the result and obtain a much
narrower gap in Section~\ref{SS:Transition} below.

\begin{figure} [htb]
\begin{picture}(52,32)
\put(0,0){\includegraphics{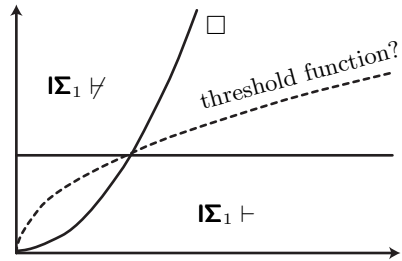}}
\put(26,30){$\square$}
\put(25,5){\small$\IS1 \vdash$}
\put(5,22){\small$\IS1 \not\vdash$}
\put(25,20){\rotatebox{15}{\small threshold function?}}
\end{picture}
\caption{\smaller\sf Transition from provability to
unprovability: $\WO\ff$ is provable from~$\IS1$ when $\ff$ is constant, and $\IS1$-unprovable when $\ff$ is the square function or above---yet it
is true; the problem of filling the gap and finding a threshold
function will be addressed in Section~\ref{SS:Transition}.}
\label{F:Threshold}
\end{figure}

\subsection{A counting formula}
\label{SS:Counting}

In order to strengthen the previous results, it will be crucial to control the
number of positive  $3$-strand braids that satisfy some constraints
simultaneously involving the complexity and the braid order. The purpose
of this section is to establish the needed estimates. Precisely, we shall count the
number of braids smaller than~$\D_3^\kk$ that have complexity at
most~$\ell$. By the results of~\cite{Dhi}, the total number of positive
$3$-strand braids with complexity at most~$\ell$ is $2^{\ell+3} - 3\ell-7$.
Discriminating according to the comparison with~$\D_3^\kk$ makes the
situation more complicated and requires a precise analysis.

\begin{prop}
\label{P:Counting}
Let $\SS_{\kk, \ell}= \{\bb \in \BB3 \mid \bb \le \D_3^\kk ~\hbox{and}~
\can\bb \le \ell\}$.
Then, for $\ell\ge\kk\ge1$, we have 
\begin{equation}
\label{E:Counting}
\card(\SS_{\kk, \ell}) = \sum_{\mm=1}^{\kk}\binom{\ell+3}{\mm+1} -
\kk + 1.
\end{equation}
\end{prop}

\begin{coro}
\label{C:Counting}
$(i)$ For all~$\ell \ge \kk \ge 1$, we have 
\begin{equation}
\label{E:Counting1}
\card(\SS_{\kk, \ell}) \le (\ell+3)^{\kk+2}.
\end{equation}

$(ii)$ For each~$\kk$, the number $\card(\SS_{\kk, \ell})$ is the value
at~$\ell$ of a degree $\kk+1$ polynomial with leading
coefficient~$1/(\kk+1)!$; in particular,  for~$\ell$ large
enough, we have 
\begin{equation}
\label{E:Counting2}
\card(\SS_{\kk, \ell}) \ge
\frac12 \ell^{\kk+1}/(\kk+1)!.
\end{equation}
\end{coro}

\begin{proof}[Proof (of Corollary~\ref{C:Counting} from
Proposition~\ref{P:Counting})]
$(i)$ The binomial $\binom{\ell+3}{\mm+1}$ is the product of
$\mm+1$~factors at most equal to~$\ell+3$, hence it is bounded above by
$(\ell+3)^{\mm+1}$, and the sum in~\eqref{E:Counting} is bounded above
by $\kk \, (\ell+3)^{\kk+1}$, hence by~$(\ell+3)^{\kk+2}$.

$(i)$ The binomial $\binom{\ell+3}{\mm+1}$ is the value at~$\ell$ of a
polynomial of degree~$\mm+1$ with leading coefficient $1/(\mm+1)!$.
By summing from $\mm = 1$ to $\mm = \kk$, we obtain a
degree~$\kk+1$ polynomial. The leading term comes from
$\binom{\ell+3}{\kk+1}$ only, so its coefficient is~$1/(\kk+1)!$.
\end{proof}

The positive $3$-strand braids with complexity at most~$\ell$ exactly
are the divisors of~$\D_3^\ell$, so the cardinality of the set~$\SS_{\kk,
\ell}$ of Proposition~\ref{P:Counting} is the rank of the braid~$\D_3^\kk$
in the $<$-increasing enumeration of the set~$\Div(\D_3^\ell)$ made by all
(left or right) divisors of~$\D_3^\ell$ in~$\BB3$. To compute this rank, we
start from the explicit description of the enumeration of~$\Div(\D_3^\ell)$
given in~\cite{Dhh}. For
$\ell\ge0$, we denote by~$\Seqi\ell$ the length $\ell+1$ sequence $(1, \ss1,
\sss12, ..., \sss1\ell)$. For $\bb$ a braid and $\Sigma$ a sequence of braids, we
write $\bb\Sigma$ for the sequence obtained by multiplying each entry
in~$\Sigma$ by~$\bb$ on the left, and, for $\Sigma, \Sigma'$ sequences of
braids, we write $\Sigma + \Sigma'$ for the concatenated sequence consisting
of~$\Sigma$ followed by~$\Sigma'$. 

\begin{lemm} 
\label{L:Enum3}
\cite[Prop.\,4.7]{Dhh}
For $\ell\ge0$, let $\dddd\ell$ denotes the braid (represented by) the length~$\ell$ suffix of
$...\sss12\sss22\sss12\ss2$, and let $\Seq\ell$ be the sequence in~$\BB3$ defined by
\begin{align} \label{E:Def3}
\Seq\ell =  \dddd0\Seqi \ell + \Seq{\ell, 1}
+ \dddd1\Seqi \ell + \cdots 
+ \dddd{2\ell-1}\Seqi \ell + \Seq{\ell, 2\ell} + \dddd{2\ell} \Seqi \ell,
\end{align}
where $\Seq{\ell, 1}, \cdots, \Seq{\ell, 2\ell}$ are 
defined by $\Seq{\ell, 1} =
\Seq{\ell, 2\ell} = \emptyset$ and, for $2 \le \mm \le 2\ell-1$, 
$$\Seq{\ell, \mm} = 
\begin{cases}
\phantom{\ss2}\ss1(\Seq{\ell-1, \mm-1} +
\dddd{\mm-1}\Seqi{\ell-1} +
\Seq{\ell-1, \mm}) & \mbox{for $\mm = 0 \pmod 4$,}\\
\ss2\ss1(\Seq{\ell-1, \mm-2} + \dddd{\mm-1}\Seqi{\ell-1} +
\Seq{\ell-1, \mm-1}) & \mbox{for $\mm = 1 \pmod 4$,}\\
\phantom{\ss2}\ss2(\Seq{\ell-1, \mm-1} +
\dddd{\mm-1}\Seqi{\ell-1} +
\Seq{\ell-1, \mm}) & \mbox{for $\mm = 2 \pmod 4$,}\\
\ss1\ss2(\Seq{\ell-1, \mm-2} + \dddd{\mm-1}\Seqi{\ell-1} +
\Seq{\ell-1, \mm-1}) & \mbox{for $\mm = 3 \pmod 4$.}
\end{cases}$$
Then $\Seq\ell$ is the $<$-increasing enumeration
of the set~$\Div(\D_3^\ell)$.
\end{lemm}

\begin{table}[t]
\begin{picture}(130,47)(0, 0)
\put(-1,10){$\dddd0 \Seqi3$}
\put(10,10){$\Seq{3,1}$}
\put(19,10){$\dddd1 \Seqi3$}
\put(30,10){$\Seq{3,2}$}
\put(39,10){$\dddd2 \Seqi3$}
\put(50,10){$\Seq{3,3}$}
\put(59,10){$\dddd3 \Seqi3$}
\put(70,10){$\Seq{3,4}$}
\put(79,10){$\dddd4 \Seqi3$}
\put(90,10){$\Seq{3,5}$}
\put(99,10){$\dddd4 \Seqi3$}
\put(110,10){$\Seq{3,6}$}
\put(119,10){$\dddd6 \Seqi3$}
\put(10,9){$\underbrace{\hbox to 26mm{\hfill}}$}
\put(13,4){$\ss2\cdot$}
\put(29,4){$\ss1\ss2\cdot$}
\put(17,3){$\swarrow$}
\put(26,3){$\searrow$}
\put(50,9){$\underbrace{\hbox to 26mm{\hfill}}$}
\put(53,4){$\ss1\cdot$}
\put(69,4){$\ss2\ss1\cdot$}
\put(57,3){$\swarrow$}
\put(66,3){$\searrow$}
\put(90,9){$\underbrace{\hbox to 26mm{\hfill}}$}
\put(93,4){$\ss2\cdot$}
\put(109,4){$\ss1\ss2\cdot$}
\put(97,3){$\swarrow$}
\put(106,3){$\searrow$}
\put(12,0){$\cdots$}
\put(31,0){$\cdots$}
\put(52,0){$\cdots$}
\put(71,0){$\cdots$}
\put(92,0){$\cdots$}
\put(111,0){$\cdots$}

\put(19,22){$\dddd0 \Seqi2$}
\put(30,22){$\Seq{2,1}$}
\put(39,22){$\dddd1 \Seqi2$}
\put(50,22){$\Seq{2,2}$}
\put(59,22){$\dddd2 \Seqi2$}
\put(70,22){$\Seq{2,3}$}
\put(79,22){$\dddd3 \Seqi2$}
\put(90,22){$\Seq{2,4}$}
\put(99,22){$\dddd4 \Seqi2$}
\put(30,21){$\underbrace{\hbox to
26mm{\hfill}}$}
\put(33,16){$\ss2\cdot$}
\put(49,16){$\ss1\ss2\cdot$}
\put(37,14.5){$\swarrow$}
\put(46,14.5){$\searrow$}
\put(70,21){$\underbrace{\hbox to
26mm{\hfill}}$}
\put(73,16){$\ss1\cdot$}
\put(89,16){$\ss2\ss1\cdot$}
\put(77,14.5){$\swarrow$}
\put(86,14.5){$\searrow$}

\put(39,34){$\dddd0 \Seqi1$}
\put(50,34){$\Seq{1,1}$}
\put(59,34){$\dddd1 \Seqi1$}
\put(70,34){$\Seq{1,2}$}
\put(79,34){$\dddd2 \Seqi1$}
\put(50,33){$\underbrace{\hbox to
26mm{\hfill}}$}
\put(53,28){$\ss2\cdot$}
\put(69,28){$\ss1\ss2\cdot$}
\put(57,26.5){$\swarrow$}
\put(66,26.5){$\searrow$}

\put(59,44){$\dddd0 \Seqi0$}
\end{picture}
\smallskip
\caption{\smaller\sf Inductive construction of~$\Seq\ell$ as a
Pascal triangle: the subsequence~$\Seq{\ell,\mm}$ is obtained by
concatenating translated copies of the previous
subsequences~$\Seq{\ell-1,\mm-1}$ and~$\Seq{\ell-1,\mm}$,
or~$\Seq{\ell-1,\mm-2}$ and~$\Seq{\ell-1,\mm-1}$, depending on
the parity of~$\mm$.}
\label{T:Triangle}
\end{table}

The result is illustrated in  Table~\ref{T:Triangle}: the
sequence~$\Seq\ell$ is constructed by starting with $2\ell+1$~copies
of~$\Seqi\ell$ translated by~$\dddd0$, ...,
$\dddd{2\ell}$ and inserting (translated copies
of) fragments of the previous sequence~$\Seq{\ell-1}$.  
For instance, we find
$\Seq0 = \dddd0 \Seqi0 =\nobreak (1)$,
which corresponds  to the trivial fact that $1$ is the only
divisor of~$1$ in~$\BB3$, then

\smallskip
$\Seq1 = \dddd0 \Seqi1 + \Seq{1,1} + \dddd1 \Seqi1 +
\Seq{1,2} + \dddd2 \Seqi1$

\hspace{1cm} $= (1, \ss1) + \emptyset
+ \ss2(1, \ss1) + \emptyset 
+ \ss1\ss2(1, \ss1)
= (1, \ss1, \ss2, \ss2\ss1, \ss1\ss2, \ss1\ss2\ss1)$,

\smallskip
\noindent which is the $<$-increasing enumeration of the
6~divisors of~$\D_3$, then
\smallskip

$\Seq2  = \dddd0 \Seqi2 + \Seq{2,1} +
\dddd1 \Seqi2 + \Seq{2,2} +
\dddd2 \Seqi2 + \Seq{2,3} +
\dddd3 \Seqi2 + \Seq{2,4} +
\dddd4 \Seqi2$

\hspace{1cm} $ = (1, \ss1, \sss12) + \emptyset
+ \ss2(1, \ss1, \sss12) + \ss2(\ss2, \ss2\ss1)
+ \ss1\ss2(1, \ss1, \sss12)$

\hspace{3cm} $ + \ss1\ss2(\ss2, \ss2\ss1)
+ \sss12\ss2(1, \ss1, \sss12) + \emptyset 
+ \ss2\sss12\ss2(1, \ss1, \sss12)$

\hspace{1cm} $ = (1, \ss1, \sss12, \ss2, \ss2\ss1, \ss2\sss12,
\sss22, \sss22\ss1,
\ss1\ss2, \ss1\ss2\ss1, \ss1\ss2\sss12,
\ss1\sss22$, 

 \hspace{3cm} $\ss1\sss22\ss1, \sss12\ss2, \sss12\ss2\ss1,
\sss12\ss2\sss12, \ss2\sss12\ss2,
\ss2\sss12\ss2\ss1, \ss2\sss12\ss2\sss12)$,

\smallskip
\noindent the $<$-increasing enumeration of the 19~divisors
of~$\D_3^2$, {\it etc.}

With the previous precise description at hand, we can now evaluate the
rank of~$\D_3^\kk$ in the sequence~$\Seq\ell$.

\begin{proof}[Proof of Proposition~\ref{P:Counting}]
With the notation of Lemma~\ref{L:Enum3}, for $1\le\mm\le\ell$, define
$\Seqt\ell\mm = \Seq{\ell,2\mm+1}+\dddd{2\mm+1}\Seqi{\ell} + 
\Seq{\ell,2\mm+2}$ (the underbraced sets in
Table~\ref{T:Triangle}). Then, by Lemma~\ref{L:Enum3}, the
sets~$\Seqt\ell\mm$ obey the inductive rules
$\Seqt11=\dddd1\Seqi1$, and, putting $\Seqt\ell\mm = \emptyset$
for $\mm\le0$ and $\mm>\ell$,
\begin{equation}
\Seqt\ell\mm = \ss{\cl\mm}\Seqt{\ell-1}{\mm-1} +
\dddd{2\mm-2}\Seqi\ell +
 \ss{\cl{\mm+1}} \ss{\cl\mm}\Seqt{\ell-1}\mm,
\end{equation}
where we recall that $\ss{\cl\mm}$ means~$\ss1$ for
$\mm$~even, and~$\ss2$ for $\mm$~odd. Left translations
in~$\BB3$ are injective, and the sequence~$\Seqi\ell$ has length
$\ell+1$, so we deduce that the length~$\seq\ell\mm$
of~$\Seqt\ell\mm$ obeys the induction rules $\seq\ell\mm=0$
for $\mm\le0$ and $\mm>\ell$, $\seq11=2$, and
$\seq\ell\mm = \seq{\ell-1}{\mm-1} + \seq{\ell-1}\mm +
\ell+1$. It follows that $\seq\ell\mm+\ell+3$ obeys the standard
Pascal triangle induction rule, and one finally obtains
\begin{equation}
\label{E:C}
\seq\ell\mm=\binom{\ell+3}{\mm+1}-\ell-3.
\end{equation}
In terms of the sequences~$\Seqt\ell\mm$, the $<$-increasing
enumeration of~$\Div(\D_3^\ell)$ is
\begin{equation}
\label{E:EnumBis}
\dddd0\Seqi\ell + \Seqt\ell1 + 
\dddd2\Seqi\ell + \Seqt\ell2 + 
\dddd4\Seqi\ell + \cdots + 
\dddd{2\ell-2}\Seqi\ell + \Seqt\ell\ell + \dddd{2\ell}\Seqi\ell.
\end{equation} 
By construction, $\D_3^\kk$ is the last element of the
sequence~$\Seq\kk$, \ie, it is~$\dddd{2\kk}\sss1\kk$---note that, by definition, we have
$\ddd\kk = \dddd{2\kk}$ for each~$\kk$, where $\ddd\kk$ is as in Lemma~\ref{L:Order}.
Now, for $\ell\ge\kk$, the element~$\dddd{2\kk}\sss1\kk$
appears in~\eqref{E:EnumBis} as the $(\kk+1)$st element of the factor $\dddd{2\kk}\Seqi\ell$,
so the rank of~$\D_3^\kk$ in~$\Seq\ell$ is the number of elements
of~\eqref{E:EnumBis} on the left of or
equal to the entry~$\dddd{2\kk}\sss1\kk$, which is
$$(\ell+1) + \seq\ell1 + (\ell+1) + \seq\ell2 + (\ell+1) +
\cdots + (\ell+1) + \seq\ell\kk + (\kk+1),$$
whence~\eqref{E:Counting} by substituting the values
given in~\eqref{E:C}.
\end{proof}

\subsection{Phase transition}
\label{SS:Transition}

We are ready for analyzing the transition between $\IS1$-provability and
$\IS1$-unprovability for~$\WO\ff$ more precisely. We start with the
provability side. We recall that $\Ack_\rr$ denotes the $\rr$th level
function in the  Grzegorczyk hierarchy, and $\Ack$ denotes the Ackermann
function, which is the diagonal function $\xx \mapsto \Ack_\xx(\xx)$. In
the sequel, we need the functional inverses of these functions: for  $\ff$ a
non-decreasing unbounded function from~$\Nat$ to itself, $\ff\inv$
denotes the function that maps~$\xx$ to the unique~$\yy$ satisfying
$\ff(\yy-1) < \xx \le \ff(\yy)$. Thus, if $\ff$ is a fast growing function, then
$\ff\inv$ is a slow growing function.

\begin{thrm}
\label{T:LowerBound}
For~$\rr \ge 0$, let $\ff_\rr$ be defined by $\ff_\rr(\xx) = 
\lfloor\!\!\sqrt[\Ack_\rr\inv(\xx)]{\xx}\rfloor$. Then the principle
$\WO{\ff_\rr}$ is provable from~$\IS1$.
\end{thrm}

\begin{proof}
As for the proof of Proposition~\ref{P:Constant}, we use a counting
argument. Let $\kk$ be a fixed number. Define $\mm = 2 \Ack_\rr(2\kk +
6)$, which makes sense inside~$\IS1$ as each function~$\Ack_\rr$ is
primitive recursive. Then we claim that $\mm$ is large enough for the result of
the principle~$\WO\ff$ to hold. 

Let $\SS = \{\bb \in \BB3 \mid  \bb \le \D_3^\kk ~\text{and}~\can\bb \le
\kk + \!\! \sqrt[2\kk+6]{\mm}\}$. With the notation of
Proposition~\ref{P:Counting}, the set~$\SS$ is the set
$\SS_{\kk, \ell}$ with $\ell = \kk + \!\! \sqrt[2\kk+6]{\mm}$, and
\eqref{E:Counting1} gives, assuming $\mm \ge 4$,
\begin{equation}
\label{E:CardS}
\card(\SS) \le (\kk + \!\! \sqrt[2\kk+6]{\mm} + 3)^{\kk+2}
\le\mm^{\frac{\kk+3}{2\kk+6}}
= \sqrt{\mm}
< \mm/2.
\end{equation}

Now, assume that $(\bb_0, ..., \bb_{\mm'})$ is a descending sequence of
braids that is $(\kk,\!\ff_\rr)$-simple. First, by hypothesis, we have $\can{\bb_0} \le \kk$, \ie,
$\bb_0$ is a divisor of~$\D_3^\kk$. As the braid order on~$\BB\infty$ extends both the left and the right
divisibility partial orders~\cite{Lve}, we deduce $\bb_0 \le \D_3^\kk$. As the sequence
$(\bb_0, ...,
\bb_{\mm'})$ is descending by hypothesis, we deduce that
$\bb_\ttii \le
\D_3^\kk$ holds for each~$\ttii$. On the other hand, consider the entries~$\bb_\ttii$ with $\ttii
\ge \mm/2$, if any. Then, by the choice of~$\mm$, we have
$$\Ack\inv_\rr(\ttii) \ge \Ack\inv_\rr(\mm/2) = 2\kk+6,$$
hence, for $\mm \ge \ttii \ge \mm/2$,  we find $\can{\bb_\ttii} \le \kk +
\!\! \sqrt[2\kk+6]{\mm}$. Thus, every such entry lies in the set~$\SS$
considered above. By \eqref{E:CardS} and the pigeonhole principle, there exist
strictly less than $\mm/2$ such braids, and, finally, we must have $\mm' \le
\mm$. 

The previous argument takes place entirely inside~$\IS1$, so we conclude that
$\WO{\ff_\rr}$ can be proved from~$\IS1$.
\end{proof}

In view of the specific form of the functions~$\ff_\rr$ involved in
Theorem~\ref{T:LowerBound}, the next  natural function to be looked at is the
one involving the inverse of the Ackermann function~$\Ack$ instead of the
functions~$\Ack_\rr$. Here comes the negative result.

\begin{thrm}
\label{T:UpperBound}
Let $\ff_\om$ be defined by $\ff_\om(\xx) =   \lfloor\!\!\sqrt[\Ack\inv(\xx)]{\xx}\rfloor$.
Then  the principle $\WO{\ff_\om}$ is not provable from~$\IS1$.
\end{thrm}

Essentially, what we do consists in replacing a constant function with the inverse of the Ackermann function.
What makes this possible is that these two functions cannot be distinguished inside~$\IS1$.
The general idea of the proof, which is reminiscent of the analysis of phase 
transition for the Kruskal theorem~\cite{Wei1, Wei2, Wei3, Wei4}, consists in
starting with a long descending sequence, typically a $\G_3$-sequence---or,
equivalently, any sequence witnessing for the principle~$\WO\square$---and
then constructing a new sequence by dilating the original one so as to lower the
complexity of the entries. The argument requires that sufficiently many braids
of low complexity are available, and this is where the estimate of
Corollary~\ref{C:Counting} is crucial.

\begin{proof}
For $\xx$ a positive integer, we write~$\logb\xx$
for $\lfloor\log_2\xx\rfloor{+}1$, \ie, for the length of the binary expansion
of~$\xx$, and we put $\logb0=0$ to complete the definition. Then we fix a
function~$\hh$, provably total in~$\IS1$---actually simply exponential---, such
that, for each~$\kk$, we have $\hh(\kk) \ge 4\kk + 10$, and $\ttii \ge
\hh(\kk)$ implies
\begin{equation}
\label{E:ChoiceH}
5\kk + 11 + (\logb\ttii)^2 +  3(\kk+1) \! \sqrt[\kk+1]{2^{\logb\ttii}}
\le \sqrt[\kk]\ttii.
\end{equation}

Let $\kk$ be a positive integer that is  large enough. Let $\bb = \D_3^\kk$,
and let
$\bb_0, \bb_1$, ... be the
$\G_3$-sequence from~$\bb$. We saw in the proof of
Theorem~\ref{T:UnprovGame} that the length of this $\G_3$-sequence is
at least~$\Ack(\kk)$, and, therefore, $\bb_\tt$ is defined for
$0 \le\nobreak \tt \le\nobreak \Ack(\kk)$. Moreover, the complexity of~$\bb$
is~$\kk$, so, using Lemma~\ref{L:Complexity} as in the proof of
Theorem~\ref{T:Jump}$(ii)$, we deduce
\begin{equation}
\label{E:Bound0}
\can{\bb_\ttii} \le \kk + 6 + \ttii^2.
\end{equation}
So, we have a descending sequence $(\bb_0, ..., \bb_{\Ack(\kk)})$ in~$\BB3$
that satisfies the complexity requirement~\eqref{E:Bound0} for~$\ttii \le
{\Ack(\kk)}$. We shall now construct a new descending sequence $(\bb'_0, ...,
\bb'_{\Ack(\kk)})$ satisfying the (much) stronger complexity requirement
\begin{equation}
\label{E:Bound}
\can{\bb'_\ttii} \le 2\hh(\kk) +
\sqrt[\Ack\inv(\ttii)]\ttii
\end{equation}
for each~$\ttii \le \Ack(\kk)$. To this end, we start from the sequence
$(\bb_{\logb0}, ... \bb_{\logb{(\Ack(\kk))}})$. 
This sequence is non-increasing---but certainly not strictly decreasing as most entries are repeated many
times. As for complexity, it is essentially $(k+6, \log)$-simple. Now, the combinatorial result of
Section~\ref{SS:Counting} will enable us to find sufficiently many braids of
low complexity which, when conveniently appended to the entries of the
previous sequence, guarantee that the final sequence is descending and
keeps the expected complexity.

Let $\SS_\ttii = \{\bb \in \BB3 \mid \bb \le \D_3^\kk ~\text{and}~
\can\bb \le (\kk+1) \sqrt[\kk+1]{2^{\logb\ttii}}  \}$. 
With the notation of Proposition~\ref{P:Counting}, the set~$\SS_\ttii$ is
$\SS_{\kk, \ell}$ with $\ell =  (\kk+1)\sqrt[\kk+1]{2^{\logb\ttii}}$, and
\eqref{E:Counting2} gives, provided $\kk$ is large enough,
$$\card(\SS_\ttii) \ge 
\frac12
\frac{(\kk+1)^{\kk+1}(\sqrt[\kk+1]{2^{\logb\ttii}})^{\kk+1}}{(\kk+1)!}
\ge 2^{\logb\ttii} \ge 2^{\logb\ttii}-\ttii.$$
Hence, for each~$\ttii$, the $<$-increasing enumeration of~$\SS_\ttii$ is a
sequence of length at least~$2^{\logb\ttii}\!-\!\ttii$ and, in particular, its
$(2^{\logb\ttii}\!-\!\ttii)$th entry (counting from~$1$) is well defined. 

We are ready to define our new sequence, \ie, to define~$\bb'_\ttii$ for $\ttii
\le \Ack(\kk)$. There are two cases. If $\ttii$ is small, namely for
$\ttii \le \hh(\kk)$, we define~$\bb'_\ttii$ to be the $3$-braid with
exponent sequence
\begin{equation}
(\underbrace{2, 2, ..., 2}_{2\kk+4\text{~entries}}, \hh(\kk)+2 - \ttii).
\end{equation}
Otherwise, \ie, for $\ttii > \hh(\kk)$, we define~$\bb'_\ttii$ to be the $3$-braid with
exponent sequence
\begin{equation}
(\ee_\pp, ..., \ee_3, \ee_2 + 1, \ee_1 + 2, 
\underbrace{2, 2, ..., 2}_{\kk+2-\qq\text{~entries}}, 
\ee'_\qq+2, \ee'_{\qq-1}, ..., \ee'_1),
\end{equation}
where $(\ee_\pp, ..., \ee_1)$ is the exponent sequence of~$\bb_{\logb\ttii}$,
and  $(\ee'_\qq, ..., \ee'_1)$ is that of the
$(2^{\logb\ttii}\!-\!\ttii)$th entry in the $<$-increasing enumeration
of~$\SS_\ttii$, which exists as observed above. The factors ``$+1$'' and
``$+2$'' are added to guarantee that the considered sequences satisfy the
normality conditions of Definition~\ref{D:Normal3}, and that the value of the
parameter~$\qq$ remains discernible. Note that the quantity $\kk+2-\qq$ is 
always nonnegative because,  by hypothesis, $\bb_{\logb\ttii} < \D_3^\kk$
holds and, therefore, by  Lemma~\ref{L:Order}, the breadth
of~$\bb_{\logb\ttii}$ is at most~$\kk+2$.

We claim that the sequence $(\bb'_0, ..., \bb'_{\Ack(\kk)})$ has the
expected properties. First, it is descending. Indeed,  for $\tti < \ttii \le
\hh(\kk)$, Proposition~\ref{P:Order3} implies $\bb'_\tti > \bb'_\ttii$
because $\bb'_\tti$ and $\bb'_\ttii$ have the same breadth and the same first
$2\kk+4$~exponents, while the last entry in the exponent sequence
of~$\bb'_\tti$ is larger than that of~$\bb'_\ttii$. 

Then, for $\tti \le \hh(\kk) < \ttii$, Proposition~\ref{P:Order3} again
implies $\bb'_\tti > \bb'_\ttii$ because the breadth of~$\bb'_\tti$, namely
$2\kk+5$, is larger than that of~Ê$\bb'_\ttii$, which is~$\pp + \kk + 2$,
hence at most $2\kk + 4$ since, as already observed above, the breadth
of~$\bb_{\logb\ttii}$ is at most~$\kk+2$. 

Next, assume $\hh(\kk) \le \tti < \ttii$ with $\logb\tti < \logb\ttii$. Then, by
hypothesis, the exponent sequence of~$\bb_{\logb\tti}$ is
$\ShortLex$-larger than that of~$\bb_{\logb\ttii}$, so
Proposition~\ref{P:Order3} implies $\bb_{\logb\tti} > \bb_{\logb\ttii}$. By
definition of the $\ShortLex$-ordering, appending
$\kk+2$ new entries at the right of the previous sequences does not change
the ordering, and, again by Proposition~\ref{P:Order3}, we deduce
$\bb'_\tti > \bb'_\ttii$. 

Finally, assume $\hh(\kk) \le \tti < \ttii$ with $\logb\tti = \logb\ttii$. Then, by
construction, the sets~$\SS_\tti$ and~$\SS_\ttii$ coincide, hence so do their
increasing enumerations. Then $\tti <\ttii$ implies
$2^{\logb\tti}-\tti > 2^{\logb\ttii} - \ttii$, and therefore again $\bb'_\tti >
\bb'_\ttii$: the result is clear if the breadth of the $(2^{\logb\tti}-\tti)$th
and $(2^{\logb\tti}-\ttii)$th entries of~$\SS_\ttii$ are equal;
otherwise, the $+2$~factor inserted in~$\ee'_\qq$ guarantees that
$\bb'_\tti > \bb'_\ttii$ holds as well. 

It remains to bound the complexity of the braids~$\bb'_\ttii$, and, for
this, it will be sufficient to use the rough connections
of~\eqref{E:Comparison}. For
$\ttii <
\hh(\kk)$, the definition of~$\bb'_\ttii$ and the hypotheses on the
function~$\hh$ give
$$\can{\bb'_\ttii} \le \lg{\bb'_\ttii} =
\hh(\kk) + 4\kk + 10 -
\ttii \le 2 \hh(\kk).$$
Assume now $\ttii \ge \hh(\kk)$. By definition, we have $\can\bb \le
(\kk+1)  \sqrt[\kk+1]{2^{\logb\ttii}}$ for each~$\bb$ in~$\SS_\ttii$, so,
by~\eqref{E:Comparison}, we deduce $\lg\bb \le 3(\kk+1) 
\sqrt[\kk+1]{2^{\logb\ttii}}$ for each such~$\bb$, whence
$$\lg{\bb'_\ttii} \le \lg{\bb_{\logb\ttii}} + 1 + 2 + 2(\kk + 2) + 2  + 
3(\kk+1)  \sqrt[\kk+1]{2^{\logb\ttii}}.$$
By construction, we always have $\lg{\bb\OP\tt} \le \lg\bb + \tt-1$, so,
iterating, we deduce
$$\lg{\bb_\tt} \le \lg{\bb_0} + 0 + 1 + \cdots + (\tt-1) < 3\kk + \tt^2.$$
Applying this to~$\bb_{\logb\ttii}$, we find 
\begin{multline*}
\qquad\lg{\bb'_\ttii}
\le 3\kk + (\logb\ttii)^2 + 2\kk + 11 + 
3(\kk+1)  \sqrt[\kk+1]{2^{\logb\ttii}} \\
= 5\kk + 11 + (\logb\ttii)^2 +  3(\kk+1) \!
\sqrt[\kk+1]{2^{\logb\ttii}}.\qquad
\end{multline*}
As $\hh$ has been chosen so as to satisfy~\eqref{E:ChoiceH}, we deduce
$\lg{\bb'_\ttii}
\le
\sqrt[\kk]\ttii$, hence $\can{\bb'_\ttii} \le \sqrt[\kk]\ttii$. Now, for $\ttii \le
\Ack(\kk)$, we have $\Ack\inv(\ttii) \le \kk$, and we finally deduce
$\can{\bb'_\ttii} \le \sqrt[\Ack\inv(\ttii)]\ttii$.
Summarizing, we conclude that, in all cases, namely $\ttii \le \hh(\kk)$ and
$\ttii>\hh(\kk)$, we have $\can{\bb'_\ttii} \le 2\hh(\kk) + 
\sqrt[\Ack\inv(\ttii)]\ttii$, \ie, \eqref{E:Bound} holds, as expected.

It is now easy to conclude. Indeed, assume that $\WO{\ff_\om}$ is
provable from~$\IS1$. This implies that there exists a function~$\gg$,
provably total in~$\IS1$ and, therefore, primitive recursive, such that, for
each~$\kk$, each $(\kk,\!\ff_\om)$-simple descending sequence has length
at most~$\gg(\kk)$. But we showed above that $\gg(2\hh(\kk))$ is larger
than~$\Ack(\kk)$ for all~$\kk$. This is impossible, as $\hh$ is
primitive recursive, the composition of two primitive recursive functions
is primitive recursive, and the Ackermann function cannot be bounded
above by any primitive recursive function. Hence
$\WO{\ff_\om}$ is not provable from~$\IS1$.
\end{proof}

\section{Extension to arbitrary braids}
\label{S:Extension}

So far, we considered $3$-strand braids and the well-order on~$\BB3$; as the latter has ordinal type~$\om^\om$, we naturally found connection with the Ackermann function and the system~$\IS1$. We shall now discuss the extension of the previous approach to arbitrary braids in~$\BBi$. As the well-order on~$\BBi$ 
has ordinal type~$\om^{\om^\om}$, we shall jump to the next level in the approximations to
the Peano system, namely the system~$\IS2$ where the induction scheme is asserted for all
$\Sigma^0_2$~sentences. The main result is that we can define a convenient notion of
$\Gsp_\infty$-sequence in~$\BBi$ so that every $\Gsp_\infty$-sequence is finite, but the latter fact
cannot be established from the axioms of~$\IS2$. 

\subsection{Special braids}
\label{S:Special}

Extending the results of Section~\ref{S:Sequence3} to arbitrary positive braids 
is both easy and non-easy. The principle is easy: in order to define
$\G_\infty$-sequences, what we need is an elementary operation
$\bb \mapsto \bb\OP\tt$ that satisfies $\bb > \bb\OP\tt$---in order to
guarantee that iterated $\G_\infty$-sequences be finite---and some formula
similar to~\eqref{E:HardyConnection}---in order to allow comparison with the
Hardy hierarchy of fast growing functions on~$\Nat$. The difficulty is
that, in order to define a convenient ordinal assignment, we need a precise
control of the rank of a braid in the well-ordering of~$\BBi$. The normal form
developed by Burckel in~\cite{Bur} can be used for this purpose, but,
contrary to the case of~$\BB3$, no explicit formula is known for the rank of a
general braid in~$\BBi$. To overcome the problem, the natural solution
consists in renouncing to define $\G_\infty$-sequences starting from
arbitrary braids in~$\BBi$, but instead restricting to specific initial braids. By defining the latter
in a convenient way---and at the expense of losing generality---we shall obtain quite simple and
satisfactory proofs.  Several solutions exist. Here, we shall develop a construction that is
simple and natural, but uses in an essential way an induction on the braid
index. Let us mention the alternative construction of~\cite{Car2}: at the
expense of using a combinatorially more intricate construction based on the
Burckel's normal form, one can directly define long descending sequences
in~$\BB\nn$ without resorting to an induction on~$\nn$.

The first step in our current approach consists in defining the notion of a
special $\nn$-braid. As mentioned above, the construction uses induction
on~$\nn$, starting with the trivial case of~$\BB2$, which under the
correspondence $\ee \mapsto \sss1\ee$  is a copy of~$\Nat$. The principle is
that a special $\nn$-braid is a certain natural composition of
special $(\nno)$-braids. In the sequel, the flip automorphism~$\flip\nn$ of
the monoid~$\BB\nn$---and of the group~$B_\nn$---plays a prominent
role, as did $\flip3$ in the case of~$\BB3$.

\begin{figure} [htb]
\begin{picture}(106, 22)
\put(4,0){\includegraphics{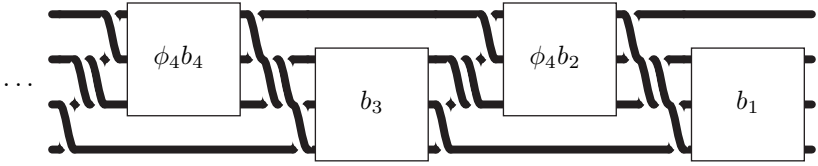}}
\put(-2,10){$\dots$}
\put(95.5,7){$\bb_1$}
\put(68,13){$\flip4\bb_2$}
\put(45.5,7){$\bb_3$}
\put(18,13){$\flip4\bb_4$}
\end{picture}
\caption{\smaller\sf The skew product of four $3$-braids~$\bb_4, ..., \bb_1$ is
the $4$-braid obtained by multiplying them after $4$-flipping each other
entry---\ie, taking the image in a horizontal medial mirror---and inserting separating patterns on
each side of flipped entries.}
\label{F:Skew}
\end{figure}

\begin{defi}
\label{D:Special}
(Figure~\ref{F:Skew}) We denote by $\flip\nn$ the {\it flip automorphism}
of~$\BB\nn$ that maps~$\ss\ii$ to~$\ss{\nn-\ii}$ for each~$\ii$.
 For $\bb_\pp, ..., \bb_1$ in~$\BB\nno$, we define
the {\it skew product} of~$\bb_\pp, ..., \bb_1$ by
\begin{equation}
\label{E:Composition}
\CC\nn\pp{\bb_\pp, ..., \bb_1} = 
\begin{cases}
\qquad \flipt\nn\bb_\pp \cdot \bb_{\pp-1} \cdot ... \cdot \flipt\nn\bb_2 \cdot
\bb_1 &\text{if $\pp$ is even,}\\
\bb_\pp \cdot \flipt\nn\bb_{\pp-1} \cdot \bb_{\pp-2} \cdot ... \cdot
\flipt\nn\bb_2
\cdot
\bb_1 &\text{if $\pp$ is odd,}
\end{cases}
\end{equation} 
with $\flipt\nn\bb = \ss1\sss22...\sss\nnd2\ss\nno \cdot \flip\nn\bb \cdot 
\ss\nno\sss\nnd2.... \sss22 \ss1$.
\end{defi}

We can now define special braids easily.

\begin{defi}
\label{D:Special}
For $\nn \ge 2$, we define an {\it $\nn$-special} braid to be, for $\nn = 2$, an
arbitrary $2$-braid, and, for $\nn \ge 3$, either the trivial braid~$1$
or a braid of the form $\CC\nn\pp{\bb_\pp, ..., \bb_1}$ where $\bb_\pp, ...,
\bb_1$ are $(\nno)$-special braids and $\bb_\pp$ is not trivial.
\end{defi}

By construction, every $(\nno$)-special braid is $\nn$-special, and,
conversely, an $\nn$-special braid belongs to~$\BB\nno$ if and only if it is
$(\nno$)-special. Hence we can drop the parameter~$\nn$ without introducing
any ambiguity and simply speak of special braids from now on.

\begin{exam}
By construction, special $3$-braids are~$1$ and those braids of the form
\begin{gather*}
\ss1 \sss2{\ee_\pp+2} \sss1{\ee_{\pp-1}+2} \, ... \,  \sss2{\ee_2+2} \sss1{\ee_1+1}
\text{\qquad with $\pp$ even and $\ee_\pp \ge 1$, $\ee_{\pp-1}, ..., \ee_1 
\ge 0$, and}\\
\sss1{\ee_\pp+1} \sss2{\ee_{\pp-1}+2} \, ... \,  \sss2{\ee_2+2} \sss1{\ee_1+1}
\text{\qquad with $\pp$ odd and $\ee_\pp \ge 1$, $\ee_{\pp-1}, ..., \ee_1 \ge 0$}.
\end{gather*}
Observe that the above expressions are $\flip{}$-normal in the sense of
Definition~\ref{D:Normal3}, that they begin and  finish with~$\ss1$, and that
every special $3$-braid is the skew product of a unique sequence of
(special) $2$-braids.
\end{exam}

It is easy to inductively extend the previous properties of special $3$-braids to
arbitrary special braids.

\begin{lemm}
\label{L:Special}
$(i)$ Every special braid has a unique word representative, which begins and
ends with~$\ss1$ whenever the braid is non-trivial.

$(ii)$ For $\nn \ge 3$, each special $\nn$-braid is the skew product of a unique
sequence of special $(\nno)$-braids.
\end{lemm}

\begin{proof}
$(i)$ Let us say that a braid word is repetitive if each letter, except possibly  the
first and the last one, is repeated at least twice and, moreover, each
letter~$\ss\ii$ is followed by a letter~$\ss\jj$ with $\vert\ii-\jj\vert \le 1$.
Owing to the braid relations of~\eqref{E:Pres}, a repetitive braid word is
equivalent to no word except itself, so any braid represented by a repetitive
word has a unique word representative. We shall inductively check that each
special $\nn$-braid has a repetitive word representative that begins and ends
with~$\ss1$. The result is obvious for $\nn = 2$. Assume that $\bb$ is
nontrivial and $\nn$-special with $\nn \ge 3$. By definition, we have $\bb =
\CC\nn\pp{\bb_\pp, ..., \bb_1}$ for some finite sequence $(\bb_\pp, ...,
\bb_1)$ of $(\nno)$-special braids. By induction hypothesis, each~$\bb_\kk$
is either trivial or it has a repetitive expression that begins and ends
with~$\ss1$. In this case, $\flip\nn\bb_\kk$ has a repetitive expression that
begins and ends with~$\ss\nno$, and, in any case, $\flipt\nn\bb_\kk$ has a
repetitive expression that begins and ends with~$\ss1$.
By~\eqref{E:Composition}, so does
$\CC\nn\pp{\bb_\pp, ..., \bb_1}$, as a product of repetitive words beginning
and ending with~$\ss1$ is still a repetitive word beginning and ending
with~$\ss1$.

$(ii)$ By~$(i)$, there is no need to distinguish between a special braid and the
unique word that represents it. Assume $\bb=\CC\nn\pp{\bb_\pp, ...,
\bb_1}$. Then the letters~$\ss\nno$ in~$\bb$ can come from the
factors~$\flipt\nn\bb_\kk$ only, and two letters~$\ss\nno$ come from the same
factor~$\flipt\nn\bb_\kk$ if and only if they are not separated by a
letter~$\ss1$. Hence, starting from~$\bb$, we recover the number of
factors~$\flipt\nn\bb_\kk$, and then each of them, and, therefore, we recover
each~$\bb_\kk$ with even~$\kk$. Finally, the factors~$\bb_\kk$ with
odd~$\kk$ are deduced, with no ambiguity on~$\pp$ because the leftmost
factor is assumed to be nontrivial.
\end{proof}

The key point in the sequel is the existence of a very simple connection 
between special braids and the braid ordering. The result is similar to what we
had in Proposition~\ref{P:Order3} with $3$-braids and their $\flip{}$-normal
form. 

\begin{defi}
Assume that $(\bb_\pp, ..., \bb_1)$ and $(\bb'_\qq, ..., \bb'_1)$ are
sequences of braids. We say that $(\bb_\pp, ..., \bb_1)$ is {\it
$\ShortLex$-smaller} than $(\bb'_\qq, ..., \bb'_1)$, denoted $(\bb_\pp, ...,
\bb_1) \lSL (\bb'_\qq, ..., \bb'_1)$, if we have either $\pp < \qq$, or $\pp =
\qq$ and there exists~$\rr$ satisfying $\bb_\rr < \bb'_\rr$ and $\bb_\kk =
\bb'_\kk$ for $\kk > \rr$.
\end{defi}

\begin{prop}
\label{P:OrderSpecial}
Assume that $\bb, \bb'$ are special $\nn$-braids with $\nn \ge 3$, say $\bb =
\CC\nn\pp{\bb_\pp, ..., \bb_1}$ and $\bb' = \CC\nn\qq{\bb'_\qq, ...,
\bb'_1}$. Then $\bb < \bb'$ holds if and only if $(\bb_\pp, ..., \bb_1)$ is {\it
$\ShortLex$-smaller} than $(\bb'_\qq, ..., \bb'_1)$.
\end{prop}

To prove this result, we need the notion of the $\BB{\nn-1}$-splitting of an
$\nn$-braid as defined in~\cite{Dho}. We recall that $\ss\ii$ is said to be a
right divisor of a positive braid~$\bb$ if there exists a positive braid~$\bb'$
satisfying $\bb = \bb' \ss\ii$. 

\begin{prop}
\label{P:Order}
\cite[Prop.\,3.8]{Dho}
$(i)$ For each~$\bb$ in~$\BB\nn$, there exists a unique sequence $(\bb_\pp, ..., \bb_1)$ in~$\BB{\nn-1}$, called the $\BB{\nn-1}$-splitting of~$\bb$, satisfying
\begin{equation}
\label{E:Flip}
\bb= \flip\nn^{\pp-1}\bb_\pp \cdot ... \cdot \flip\nn^2\bb_3
\cdot \flip\nn\bb_2 \cdot \bb_1
\end{equation}
such that, for each~$\kk\ge1$,
\begin{equation}
\label{E:NNormality}
\text{the only~$\ss\ii$ dividing $\flip\nn^{\pp-\kk}\bb_\pp \cdot ...\cdot
\flip\nn\bb_{\kk+1} \cdot \bb_\kk$  on the right is~$\ss1$}.
\end{equation}

$(ii)$ For $\nn \ge 3$ and $\bb, \bb'$ in~$\BB\nn$, the relation $\bb < \bb'$ holds if and only if the $\BB{\nn-1}$-splitting of~$\bb$ is $\ShortLex$-smaller than the $\BB{\nn-1}$-splitting of~$\bb'$.
\end{prop}

In the case $\nn = 3$, the entries in the $\BB2$-splitting of~$\bb$ are elements of~$\BB2$, \ie,
powers of~$\ss1$, and one easily checks that the exponent sequence of~$\bb$ is $(\ee_\pp, ...,
\ee_1)$ if and only if  the $\BB2$-splitting of~$\bb$ is $(\sss1{\ee_\pp}, ..., \sss1{\ee_1})$.
Thus Proposition~\ref{P:Order} directly extends Proposition~\ref{P:Order3}. 

\begin{proof}[Proof of Proposition~\ref{P:OrderSpecial}]
Assume $\nn \ge 3$, and let $\bb, \bb'$ be special $\nn$-braids. Assume 
$\bb = \CC\nn\pp{\bb_\pp, ..., \bb_1}$, $\bb' = \CC\nn\qq{\bb'_\qq, ...,
\bb'_1}$. In order to compare~$\bb$ and~$\bb'$ using  the criterion of
Proposition~\ref{P:Order}, we need to determine their $\BB{\nn-1}$-splittings.
Assume that $\pp$ is even. Put $\tau_\nn = \sss\nnd2 \,...\, \sss22 \ss1$.
Then, applying the definition of~$\flipt\nn$, we find
\begin{equation}
\label{E:Splitting}
\bb = (\ss1) \cdot \flip\nn(\tau_\nn \bb_\pp \ss1) \cdot (\tau_\nn \bb_{\pp-1} \ss1) \cdot \, ... \, \cdot \flip\nn(\tau_\nn \bb_2 \ss1) \cdot (\tau_\nn \bb_1).
\end{equation}
We claim that \eqref{E:Splitting} displays the $\BB{\nn-1}$-splitting of~$\bb$. Indeed, the right hand side term consists of factors that alternatively belong to~$\BB\nno$ and~$\flip\nn\BB\nno$, 
and, as the word is equivalent to no other word than itself, the divisibility
condition of~\eqref{E:NNormality} is satisfied. If $\pp$ is odd, 
\eqref{E:Splitting} is to be replaced with
\begin{equation}
\label{E:Splitting2}
\bb = (\bb_\pp\ss1) \cdot \flip\nn(\tau_\nn \bb_{\pp-1} \ss1) \cdot (\tau_\nn \bb_{\pp-2} \ss1) \cdot \, ... \, \cdot \flip\nn(\tau_\nn \bb_2 \ss1) \cdot (\tau_\nn \bb_1),
\end{equation}
and the result is similar. 

Assume for instance that both $\pp$ and~$\qq$ are even. Applying Proposition~\ref{P:Order}, we see that $\bb < \bb'$ holds if and only if we have either $\pp < \qq$ or 
\begin{multline}
\label{E:Compar}
\qquad (\ss1, \tau_\nn \bb_\pp \ss1, \tau_\nn \bb_{\pp-1} \ss1,  \, ... \, , \tau_\nn \bb_2 \ss1, \tau_\nn \bb_1) \\ \lSL (\ss1, \tau_\nn \bb'_\pp \ss1, \tau_\nn \bb'_{\pp-1} \ss1,  \, ... \, , \tau_\nn \bb'_2 \ss1, \tau_\nn \bb'_1).\qquad
\end{multline}
Now we observe that $\tau_\nn \bb_\kk \ss1 = \tau_\nn \bb'_\kk \ss1$ is equivalent to $\bb_\kk = \bb'_\kk$, because the monoid~$\BB\nn$ admits left and right cancellation, and that $\tau_\nn \bb_\kk \ss1 < \tau_\nn \bb'_\kk \ss1$ is equivalent first to $\bb_\kk \ss1 < \bb'_\kk \ss1$, because the order~$<$ is compatible with multiplication on the left, and then to $\bb_\kk < \bb'_\kk$, because $\bb \ss1$ is always the immediate successor of~$\bb$ in the braid ordering. So \eqref{E:Compar} is equivalent to $(\bb_\pp, ..., \bb_1) \lSL (\bb'_\pp, ..., \bb'_1)$, as expected.

The argument is similar if $\pp$ and $\qq$ are odd, and if they have different parities, owing to the fact that $\bb_\pp$ and $\bb'_\qq$ are not trivial.
\end{proof}

\subsection{$\Gsp_\nn$-sequences}

With the notion of a special braid at hand, we can now mimick the
approach of Section~\ref{S:Sequence3} and define long descending sequences
in~$\BB\nn$. Once again, the construction uses induction on the braid
index~$\nn$, and the principle is the same as for $\G_3$-sequences. 

\begin{defi}
\label{D:Descent}
$(i)$ For $\nn \ge 2$ and $\tt \ge1$, we define $\theta_{\nn,\tt} =
\CC\nn\tt{\ss1, 1, ..., 1}$.

$(ii)$ For $\nn \ge 2$, $\bb$ a nontrivial special
$\nn$-braid, and $\tt \ge 1$, the braid~$\bb\OPsp\tt_{\!\nn}$  is defined, for
$\nn = 2$, to be the braid obtained from~$\bb$ by removing one
letter~$\ss1$, and, for $\nn \ge 3$ and $\bb =
\CC\nn\pp{\bb_\pp, .., \bb_\rr, 1, ..., 1}$, putting $\bb'_\rr
= \bb_r\OPsp\tt_{\!\nno}$, to be the braid
$$\begin{cases}
\CC\nn\pp{\bb_\pp, ..., \bb_{\rr+1}, \bb'_\rr, 1, ... , 1}, 
&\text{for $\rr = 1$ or $\bb_\rr \not= \bb'_\rr \ss1$},\\ 
\CC\nn\pp{\bb_\pp, ..., \bb_{\rr+1}, \bb'_\rr, \theta_{\nno,\tt}, 1, ... , 1}, 
&\text{for $\rr \ge 2$ and $\bb_\rr = \bb'_\rr \ss1$ with $\bb'_\rr \not= 1$
or
$\pp > \rr$},\\ 
\CC\nn{\pp-1}{\theta_{\nno,\tt}, 1, ... , 1}, 
&\text{for $\pp = \rr \ge 2$ and $\bb_\rr = \ss1$}.
\end{cases}$$
Finally, we define the {\it $\Gsp_\nn$-sequence
from~$\bb$} to be the sequence $(\bb_0, \bb_1, ...)$ defined by $\bb_0 =
\bb$ and $\bb_\tt = \bb_{\tt-1}\OPsp\tt_{\!\nn}$; the sequence stops when the
trivial braid~$1$ is possibly obtained.
\end{defi}

The idea is simple: $\bb\OPsp\tt_{\!\nn}$ is obtained from~$\bb$ by considering
the rightmost nontrivial component~$\bb_\rr$ in the decomposition of~$\bb$
as a skew product, and replacing it with $\bb_r\OPsp\tt_{\!\nno}$, \ie, in
applying the rule inductively; now, if going from~$\bb_\rr$
to~$\bb_r\OPsp\tt_{\!\nno}$ amounts to deleting the last letter
in~$\bb_\rr$---necessarily a $\ss1$ according to Lemma~\ref{L:Special}---and if
$\rr$ is at least~$2$, then we add~$\theta_{\nn,
\tt}$ in the next component. In the particular case $\nn = 3$, we nearly recover
the rule of Section~\ref{S:Sequence3}. Indeed, going from~$\bb_\rr$
to~$\bb_r\OPsp\tt_{\!2}$ simply means removing one~$\ss1$ in~$\bb_\rr$,
and, then we have $\theta_{2,\tt} = \sss1\tt$, so adding~$\theta_{2, \tt}$
amounts to adding $\tt$~letters in the next block, whenever the latter is not the
final block of~$\ss1$'s. Thus the only difference between $\G_3$- and
$\Gsp_3$-sequences is that, in the latter, the separating patterns $\ss1\ss2$
and $\ss2\ss1$ play a specific role.

\begin{exam}
Let $\bb = \ss1\sss24\ss1$. Then $\bb$ is a special $3$-braid, corresponding
to the skew product $\CC32{\sss12, 1}$. Then $\bb\OPsp1_{\!3}$ is obtained
by applying the rule of~$\Gsp_2$ to the rightmost nontrivial component
of~$\bb$, namely~$\sss12$, hence removing one~$\ss1$ there. As the
parameter~$\rr$ of Definition~\ref{D:Descent} is~$2$ here, we add the
factor~$\theta_{2,1}$, \ie, $\ss1$ in the next component, so $\bb\OPsp1_{\!3}$ is
$\CC32{\ss1,\ss1}$, \ie,
$\ss1\sss23\sss12$. Iterating the process, we find the $\Gsp_3$-sequence
$$\ss1\sss24\ss1, \quad \ss1\sss23\sss12, \quad \ss1\sss23\ss1, \quad
\sss13, \quad \sss12, \quad \ss1, \quad 1.$$
\end{exam}

In the general case, the factor~$\theta_{\nno,\tt}$ that is added is more
complicated than just a power of some~$\ss\ii$. For instance, for $\nn = 3$,
the successive braids~$\theta_{3,\tt}$ turn out to be $\ss1$, $\ss1\sss22\ss1$,
$\sss12 \sss22 \ss1$, $\ss1 \sss23 \sss12 \sss22 \ss1$, etc. Some
flexibility exists here. The current values have been chosen so as to guarantee an
easy connection with the subsequent ordinal assignment.

Before investigating $\Gsp_\nn$-sequences more precisely, let us observe
that, if $\bb$ is a special $(\nno)$-braid, then, by definition, we
have $\bb\OPsp\tt_{\!\nn} = \bb\OPsp\tt_{\nno}$ for each~$\tt$. So, once again,
we can skip the index~$\nn$ without ambiguity. Inductively, there is no need
to distinguish between $\Gsp_\nn$- and $\Gsp_\nno$-sequences, and
we refer from now to $\Gsp_\infty$-sequences for all such sequences, in
the same way as $\BBi$ is seen as the union of all~$\BB\nn$'s.

\subsection{Finiteness of $\Gsp_\infty$-sequences}

As in the case of $\G_3$-sequences, we observe that, although very long
$\Gsp_\infty$-sequences exist, no such sequence is infinite, \ie, we establish
the counterpart to Proposition~A of the introduction.

\begin{prop}
\label{P:TerminationN}
For each special braid~$\bb$, the $\Gsp_\infty$-sequence from~$\bb$ is
finite,
\ie, there exists a finite number~$\tt$ satisfying $\bb\OPsp1\OPsp2... \OPsp\tt = 1$. 
\end{prop}

Like Proposition~\ref{P:Termination}, Proposition~\ref{P:TerminationN}
directly  follows from the conjunction of two results, namely that, according to
Theorem~\ref{T:Order}$(ii)$, the braid order on~$\BBi$ is a well-order, and
that every $\Gsp_\infty$-sequence is descending with respect to
that order. The latter is a consequence of

\begin{lemm}
\label{L:DescendingN}
For each special braid~$\bb$ in~$\BB\nn$ and each~$\tt$, we have
$\bb > \bb\OPsp\tt$.
\end{lemm}
\begin{proof}
An obvious induction on~$\nn$, applying the criterion of Proposition~\ref{P:OrderSpecial} to the explicit construction
of Definition~\ref{D:Descent}.
\end{proof}

\subsection{An unprovability result}

We turn to the counterpart of Theorem~A, and prove that
the finiteness of $\Gsp_\infty$-sequences cannot be proved in the
system~$\IS2$. 

\begin{thrm}
\label{T:UnprovGameN}
Proposition~\ref{P:TerminationN} is an arithmetic statement that cannot be 
proved from the axioms~$\IS2$.
\end{thrm}

As in the case of~$\BB3$, Theorem~\ref{T:UnprovGameN} follows from the result
that  the function measuring the length of $\Gsp_\infty$-sequences in terms
of the size of the initial braid grows faster than any function whose totality is
provable from~$\IS2$, and the proof relies on assigning convenient ordinals to
special braids. 

\begin{defi}
For $\bb$ a special $\nn$-braid, we define $\ordsp_\nn(\bb)$ by
$\ordsp_2(\sss1\ee) = \ee$, and 
\begin{equation}
\label{E:OrdinalN}
\ordsp_\nn(\bb) = \om^{\om^\nnd \cdot (\pp-1)} \cdot \ordsp_\nno(\bb_\pp) + \cdots +
\om^{\om^\nnd} \cdot \ordsp_\nno(\bb_2) + \ordsp_\nno(\bb_1)
\end{equation} 
for $\bb = \CC\nn\pp{\bb_\pp, ..., \bb_1}$.
\end{defi}

We observe that, for $\bb$ an $(\nno)$-braid, we  have $\ordsp_\nn(\bb) =
\ordsp_\nno(\bb)$, and,  therefore, dropping the $\nn$~subscripts introduces
no ambiguity. Formula~\eqref{E:OrdinalN} is more simple
than its counterpart~\eqref{E:Ordinal3} of Section~\ref{S:Sequence3} because,
by restricting to special braids, we avoid the problem of counting how many
normal words lie below a given one.

The construction of distinguished cofinal sequences given in
Definition~\ref{D:Fundamental3} easily extends to ordinals
below~$\om^{\om^\om}$---and even below~$\eps_0$. We recall that $\eqCNF$ refers to the Cantor Normal
Form.

\begin{defi}
\label{D:HardyN}
For $\l$ a limit ordinal below~$\eps_0$, we put
$$
\l[\xx]:=
\begin{cases}
\gamma + \om^{\delta}\cdot \xx &
\text{for $\l \eqCNF \g + \om^{\delta+1}$},\\
\gamma + \om^{\delta[\xx]} &
\text{for $\l \eqCNF \g + \om^\delta$ with $\delta$ a limit ordinal}.
\end{cases}
$$
\end{defi}

As previously, we extend to non-limit ordinals by $0[\xx]= 0$ and $(\alpha+1)[\xx]=\nobreak\alpha$
for every~$\xx$. Once fundamental sequences have been defined for all ordinals
below~$\eps_0$, the hierarchy of Hardy functions~$H_\a$ is introduced by
extending the defining relations~\eqref{E:Hardy3} to each ordinal
below~$\eps_0$. Then, as in Section~\ref{S:Sequence3}, we have the
following connection:

\begin{lemm}
\label{L:OrdConnectionN} 
For every nontrivial special braid $\bb$ and every $\tt$ in~$\Nat$, we have
\begin{equation}
\label{E:OrdConnectionN}
\ordsp(\bb\OPsp\tt) = \ordsp(\bb)[\tt].
\end{equation}
\end{lemm}

\begin{proof}
Everything has been done so as to guarantee the connection. As a preliminary step, we first check
that, when we go from~$\bb$ to~$\bb\OPsp\tt$, the case $\bb = \bb\OPsp\tt \cdot \ss1$ occurs if
and only if $\ordsp(\bb)$ is a successor ordinal. We use induction on~$\nn$. For $\nn = 2$, the
equivalence is obvious. Assume $\nn \ge 3$. Let $\bb'$ stand for $\bb\OPsp\tt$. Write $\bb = 
\CC\nn\pp{\bb_\pp,...,\bb_\rr,1,...,1}$ and $\bb' =  \CC\nn{\pp'}{\bb'_{\pp'},...,\bb'_1}$.
Assume first that $\ordsp(\bb)$ is a successor. We claim that $\bb=\bb' \cdot \ss1$ holds. If
$\ordsp(\bb)$ is a successor, then, by definition of the ordinal assignment, we must have $\rr=1$
and $\ordsp(\bb_\rr)$ is a successor, for, otherwise, the rightmost term in~$\ordsp(\bb)$ is at
least $\om^{\om^{\nn-2}}$. Since $\bb_1$ is a special $(\nn-1)$-braid, the induction
hypothesis implies $\bb_1 = \bb'_1 \cdot \ss1$. Then, by definition of~$\bb'$, we have
$$\bb'\ss1 = \CC\nn{\pp'}{\bb'_{\pp'},...,\bb'_1} \cdot \ss1 = ... \ \bb_1' \ss1 = \bb.$$
Conversely, assume $\bb=\bb' \cdot \ss1$. By the rules of the game, $\bb'$ is obtained from
$\bb$ by deletion of a rightmost~$\ss1$. This is the case only if we have $\bb =
\CC\nn\pp{\bb_\pp, ... , \bb_\rr,1, ... ,1}$ with  $\rr=1$. If $\ordsp(\bb_\rr)$ is a successor, we
are done. Otherwise, by induction hypothesis, we have $\bb_\rr \neq \bb_\rr'\ss1$. Since
$\bb'$ is obtained by replacing $\bb_\rr$ by $\bb_\rr'$ as the rightmost component
in the skew product defining~$\bb$, this contradicts the hypothesis $\bb =\bb'\ss1$.

We can now establish Equality~\eqref{E:OrdConnectionN}, distinguishing between the three cases
of Definition~\ref{D:Descent}, of which we adopt the notation. Assume first $\rr =1$ or $\bb_\rr
\neq \bb_\rr'\ss1$. The case $\rr = 1$ corresponds to
$\bb = \CC\nn\pp{b_p, ... ,b_1}$ and $b\OPsp{t} = \CC\nn\pp{b_p, ... ,b'_1}$,
and we find
\begin{align*}
\ordsp(b) 
& = \om^{\om^{n-2}{\cdot}(p-1)} \cdot \ordsp(b_p)+  ...  + \ordsp(b_1),\\
\ordsp(b\OPsp{t}) 
& = \om^{\om^{n-2}{\cdot}(p-1)}\cdot\ordsp(b_p)+  ...  +
\ordsp(b'_1)\\
& = \om^{\om^{n-2}{\cdot}(p-1)}\cdot\ordsp(b_p)+  ...  + \ordsp(b_1)[t]
= \ordsp(b)[t],
\end{align*}
where the second equality holds by inductive hypothesis.

Similarly, the case $r > 1$ and $b_r \neq b_r'\ss1$ corresponds to
$b= \CC\nn\pp{b_p, ... ,b_r,1, ... ,1}$ and $b\OPsp{t} = \CC\nn\pp{b_p, ... ,b'_r,1, ... ,1}$, 
and we find now
\begin{align*}
\ordsp(b) 
& = \om^{\om^{n-2}{\cdot}(p-1)} \cdot \ordsp(b_p)+  ...  +
\om^{\om^{n-2}{\cdot}(r-1)}\ordsp(b_r),\\
\ordsp(b\OPsp{t}) & = \om^{\om^{n-2}(p-1)}\cdot\ordsp(b_p)+  ...  +
\om^{\om^{n-2}{\cdot}(r-1)}\cdot\ordsp(b'_r)\\
& = \om^{\om^{n-2}{\cdot}(p-1)}\cdot\ordsp(b_p)+  ...  +
\om^{\om^{n-2}{\cdot}(r-1)}\cdot(\ordsp(b_r)[t])\\
& = \om^{\om^{n-2}{\cdot}(p-1)}\cdot\ordsp(b_p)+  ...  +
(\om^{\om^{n-2}{\cdot}(r-1)}\cdot\ordsp(b_r))[t]
= \ordsp(b)[t],
\end{align*}
where the second equality holds by inductive hypothesis, and the third one is true because, 
according to the preliminary result, $\ordsp(b'_r)$ is a limit ordinal.

Assume now $r\geq 2$ and $b_r = b_r'\ss1$ with $b_r'\neq 1$ or $p>r$. In this case, we find
$b= \CC\nn\pp{b_p, ... ,b_r,1, ... ,1}$ and $b\OPsp{t} =
\CC\nn\pp{b_p, ... ,b_r,\theta_{n-1,t+1}, ... ,1}$, 
leading to
\begin{align*}
\ordsp\!(b) 
& = \om^{\om^{n-2}{\cdot}(p-1)}\cdot\ordsp\!(b_p)+  ...  +
\om^{\om^{n-2}{\cdot}(r-1)}\cdot\ordsp\!(b_r),\\
\ordsp\!(b\OPsp{t}) & = \om^{\om^{n-2}{\cdot}(p-1)}\cdot\ordsp\!(b_p)+  ...  +
\om^{\om^{n-2}{\cdot}(r-1)}\cdot\ordsp\!(b'_r)\\
&\hspace{65mm} +\om^{\om^{n-2}{\cdot}(r-2)}\cdot\ordsp\!(\theta_{n-1,t})\\
& = \om^{\om^{n-2}{\cdot}(p-1)}\cdot\ordsp\!(b_p)+  ...  +
\om^{\om^{n-2}{\cdot}(r-1)}\cdot\ordsp\!(b'_r) \\
&\hspace{65mm} +\om^{\om^{n-2}{\cdot}(r-2)+\om^{n-3}\cdot t}\\
& = \om^{\om^{n-2}{\cdot}(p-1)}\cdot\ordsp\!(b_p)+  ...  +
\om^{\om^{n-2}{\cdot}(r-1)}\cdot(\ordsp\!(b_r)[t]) \\
&\hspace{65mm} +
\om^{\om^{n-2}{\cdot}(r-2)+\om^{n-3}\cdot t}\\
& = \om^{\om^{n-2}{\cdot}(p-1)}\cdot\ordsp\!(b_p)+  ...  +
\om^{\om^{n-2}{\cdot}(r-1)}\cdot(\ordsp\!(b_r)[t])\\
&\hspace{65mm} + \om^{\om^{n-2}{\cdot}(r-1)}[t]\\
& = \om^{\om^{n-2}{\cdot}(p-1)}\cdot\ordsp\!(b_p)+  ...  +
(\om^{\om^{n-2}{\cdot}(r-1)}\cdot\ordsp\!(b_r))[t]
= \ordsp\!(b)[t],
\end{align*}
where the third equality holds by induction hypothesis, and the penultimate one holds
because $\ordsp\!(b_r)$ is a successor ordinal, as established in the preliminary step.

Assume finally $p=r\ge 2$ and $b_r = \ss1$. We have
$b= \CC\nn r{b_r,1, ... ,1}$ and $b\OPsp{t} =
\CC\nn{r-1}{\theta_{n-1,t+1},1, ... ,1}$, 
and we find now
\begin{align*}
\ordsp(b) 
& = \om^{\om^{n-2}{\cdot}(r-1)}\cdot\ordsp(b_r)
 = \om^{\om^{n-2}{\cdot}(r-1)}\cdot\ordsp(\ss1)
 = \om^{\om^{n-2}{\cdot}(r-1)},\\
\ordsp(b\OPsp{t}) 
& = \om^{\om^{n-2}{\cdot}(r-2)}\cdot\ordsp(\theta_{n-1,t})
 = \om^{\om^{n-2}{\cdot}(r-2)+\om^{n-3}\cdot t}
 = \ordsp(b)[t].
\end{align*}
So \eqref{E:OrdConnectionN} holds in every case, as expected.
\end{proof}

We easily deduce a comparison between the length of
$\Gsp_\infty$-sequences and the functions of the Hardy hierarchy.
Let $\TTsp(\bb)$ denote the length of the $\Gsp_\infty$-sequence
from~$\bb$. Using exactly the same argument as for Proposition~\ref{P:Hardy}, we obtain:

\begin{prop}
\label{P:HardyN} 
Assume that $\bb$ is a special braid satisfying $\ordsp(\bb) = \b$. Then, for each~$\kk$, we have
\begin{equation}
\label{E:HardyN}
\TTsp(\bb\sss1{\kk}) = H_\b(\kk+1)-1.
\end{equation}
\end{prop}

Then we conclude as in Section~\ref{S:Sequence3}:

\begin{proof}[Proof of Theorem~\ref{T:UnprovGameN}]
Put $\bb_\kk = \ss1\sss22...\sss{\kk+1}2\sss{\kk+2}3\sss{\kk+1}2 ...
\sss22\ss1$, \ie, $\bb_\kk = \CC{\kk+3}2{\ss1,1}$. Then $\bb_\kk$ is the
smallest special $(\kk+3)$-braid which is not a $(\kk+2)$-braid. An easy
induction using~\eqref{E:OrdinalN} gives the equality $\ordsp(\bb_\kk) =
\om^{\om^\kk}$ for each~$\kk$. Now define $\UUsp(\kk) =
\TTsp(\bb_\kk \sss1\kk) + 1$. Then 
\eqref{E:HardyN} plus the definition of~$H_{\om^{\om^\om}}$ give 
$$\UUsp(\kk) = H_{\om^{\om^\kk}}(\kk+1) = H_{\om^{\om^\om}}(\kk).$$  
Therefore, the function~$\UUsp$ is~$H_{\om^{\om^\om}}$, a recursive
function that is not provably total in~$\IS2$. Now, if the finiteness of
$\Gsp_\infty$-sequences were provable in~$\IS2$, the  function~$\UUsp$
would be provably total in $\IS2$. 
\end{proof}

\subsection{Further questions}

Above, we extended the results of Section~\ref{S:Sequence3} about
$\G_3$-sequences involving $3$-strand braids to general
$\Gsp_\infty$-sequences involving arbitrary braids. We conjecture that the
results of Section~\ref{S:Transition} about arbitrary long descending sequences
in~$(\BB3, <)$ might be similarly extended to~$(\BBi, <)$. At the moment, the missing piece is a
combinatorial result analogous to Proposition~\ref{P:Counting} in the general
case. Even the cardinality of the set of all $\nn$-braids with complexity at
most~$\ell$ remains out of control at the moment, and the results of~\cite{Dhh,
Dhi} seem to discard any possibility of explicitly describing the $<$-increasing
enumeration of the set above. However, what is needed for the proof of
Theorems~\ref{T:LowerBound} and~\ref{T:UpperBound} are the rather rough
estimates of Corollary~\ref{C:Counting}, and it is not hopeless to establish similar bounds in the general case.

Also, one could look to phase transitions of a different type. For $\ff$ a function of~$\Nat$ into itself, we
may consider the variant of the $\G_3$-sequence in which $\ff(\tt)$ new crossings appear at Step~$\tt$,
instead of~$\tt$. For which functions~$\ff$ do we obtain provability/unprovability in~$\IS1$? The threshold
result turns out to be the same as for the principle~$\WO\ff$ of Section~\ref{S:Transition}.
Moreover, a similar threshold result can be obtained without much difficulty for the extension to~$\BBi$,
which---as was said above---is not known so far in the case of~$\WO\ff$.

Other natural questions involve alternative braid orders. There exists a
large space of linear orders on the braid groups~$B_\nn$ that are
compatible with multiplication on one side. Most of them do not
induce well-orderings on the braid monoid~$\BBi$, but but
at least all the orderings stemming from the hyperbolic geometry approach
suggested by W.\,Thurston and investigated in~\cite{ShW} do. For each of
these (uncountably many) orderings, and, in particular, for those (countably
many) for which there exists a more or less explicit description, one might
investigate the associated ordinal. It would be interesting to know whether
$\om^{\om^{\nn-2}}$ appears for each such order on~$\BB\nn$.

Similar questions arise when, instead of considering the braid
monoids~$\BB\nn$, one consider the dual monoids of~\cite{BKL}. Recent results by J.\,Fromentin suggest that
the restriction of the standard braid ordering to the $\nn$-strand dual braid monoid is a well-order of
ordinal type~$\om^{\om^{\nn-2}}$, and all results mentioned in the current
paper are likely to extend to the Birman--Ko--Lee context.

\section*{Appendix: Basic definitions from logic}

\subsection*{Ordinals}

Once the existence of at least one infinite set is assumed, the basic
properties of sets as captured in the Zermelo--Fraenkel system~$\ZF$
guarantee the existence of an infinite sequence of objects called ordinals,
which can be seen as a proper end-extension of the sequence of natural
numbers. Ordinals come equipped with a canonical well-order. The smallest
infinite  ordinal is denoted by~$\om$, so, by construction, the ordinals that are
smaller than~$\om$ are (a copy of) the natural numbers.

The ordinals are equipped with arithmetic operations, addition,
multiplication, exponentiation, that extend those of natural numbers, and
that obey natural algebraic laws, associativity, distributivity, etc.---but neither commutativity nor right
cancellativity. For each ordinal~$\a$, the ordinal~$\a+1$ is the immediate successor of~$\a$ in
the well-ordering of ordinals and, for instance, $\om + \om$, which is also~$\om
\cdot 2$, is the supremum of the ordinals~$\om + \kk$ with $\kk$ a natural
number. Similarly, $\om^2$ is the supremum of the ordinals~$\om\cdot
\kk$ for $\nn$ a natural number, and $\om^\om$---the ordinal type of the
well-order on~$\BB3$---is the supremum of the sequence $1, \om, \om^2, ...,
\om^\kk, ...$, while $\om^{\om^\om}$ ---the ordinal type of the
well-order on~$\BB\infty$---is the supremum of the sequence $1, \om,
\om^\om, \om^{\om^2}, \om^{\om^3}, ...$ Finally, one denotes by~$\eps_0$ the
supremum of the sequence $1, \om, \om^\om, \om^{\om^\om}, ...$ It should be kept
in mind that the ordinal~$\eps_0$, despite being very large, is countable, as well as
all other ordinals mentioned above.

When ordinals are equipped with the order topology, those of the form
$\a+1$, naturally called successor ordinals, are isolated points, while those not
of that form are limit points, and they are called limit ordinals. A
positive ordinal is limit if and only if it can be written as $\om \cdot \a$
for some ordinal~$\a$.

All ordinals mentioned so far---and many more---are constructive in the precise sense that their structure
(their build-up) can be described and their order relations decided recursively, and even elementary
recursively. In this way, constructive transfinite ordinals can be expressed and manipulated in first-order
arithmetical systems. For our concerns it is sufficient to know that there exists recursive, even elementary
recursive, ordinal notation systems for ordinals below~$\eps_0$. The standard such system is based on the
idea of the Cantor Normal Form, which we now describe.

Every ordinal~$\a$ below $\eps_0$ has a unique expression of the form
$\om^{\a_\pp}\cdot \kk_\pp +  ...  + \om^{\a_2}\cdot \kk_2 +
\om^{\a_1}\cdot \kk_1$ with  $\a > \a_\pp >  ...  > \a_2 > \a_1  \ge
0$ and $\kk_\pp, ... ,\kk_1 > 0$. This expression is called the Cantor Normal
Form of~$\a$. For $\a < \om^\om$, all exponents in the Cantor
Normal Form of~$\a$ are natural numbers, and, for $\a < \om^{\om^\om}$,
they are ordinals below~$\om^\om$.  

\subsection*{Peano Arithmetic and its subsystems}

The standard axiomatic system for formalizing arithmetic is the Peano
system~$\PA$. It deals with so-called first order arithmetic sentences,
which involve two constants~$0, 1$, two binary operations~$+, *$, and a binary relation~$<$. The
term ``first order'' means that we restrict to formulas in which all variables
refer to integers---typically, no variable may refer to a set of integers.
Then $\PA$ consists of the axioms
\begin{gather*}
1\not=0,\quad
\forall\xx(\xx+1\not=0), \quad
\forall\xx,\yy(\xx+1=\yy+1\Rightarrow\xx=\yy)\\
\forall\xx(\xx+0=\xx)\quad
\forall\xx,\yy(\xx+(\yy+1)=(\xx+\yy)+1)\\
\forall\xx(\xx*0=0)\quad
\forall\xx,\yy(\xx*(\yy+1)=(\xx*\yy)+\xx)\\
\forall\xx,\yy(\xx < \yy \Leftrightarrow \exists\zz\not= 0(\yy = \zz + \xx)),
\end{gather*}
here denoted~$\PA_0$, plus the induction axiom
\begin{equation}
\label{E:Ind}
\tag{$\Ind_{\For}$}
\forall\yy_1,...,\yy_\pp((\For(0)\ \&\ \forall\xx(\For(\xx)\Rightarrow\For(\xx+1)))
\ \Rightarrow\ 
\forall\xx(\For(\xx)))
\end{equation}
for each first order arithmetical formula~$\For$ with free variables among $\xx,
\yy_1,...,\yy_\pp$.

For~$\kk$ a natural number, we say that a formula~$\For$ is a $\Sigma^0_\kk$ if it is
equivalent to some formula
$$\exists\xx_1
\forall\xx_2 \exists\xx_3 ... Q\xx_\kk(\Forr(\xx_1, ..., \xx_\kk, \xx))$$
where all quantifiers in~$\Forr$ are bounded quantifiers, \ie, are of the form $\forall\xx < \yy$ and
$\exists\xx < \yy$.  Then $\IS\kk$ denotes the subsystem of~$\PA$
consisting of the axioms~$\PA_0$, plus the induction
axiom~$\Ind_{\For}$ for each $\Sigma^0_\kk$ formula~$\For$. So, by definition, we have
$\PA = \bigcup_{\kk\in\Nat}\IS\kk$, but, for each fixed~$\kk$, the system~$\IS\kk$, containing less axioms
than~$\PA$, is less powerful than~$\PA$: {\it a priori}, less sentences can be proved from the axioms
of~$\IS\kk$ than form the axioms of~$\PA$. It is known that, for each~$\kk$, the inclusion is in fact proper.

We say that a function~$\ff : \Nat\to \Nat$ is provably total in a formal
system~$\SS$ if there is a $\Sigma^0_1$-formula~$\For$ such
that  $\yy = \ff(\xx)$ is equivalent to~$\For(\xx, \yy)$ and there is a formal
proof of the sentence $\forall\xx\exists\yy(\For(\xx,\yy))$ from the axioms
of~$\SS$. There is a close connection between the logical strength of a formal
system~$\SS$ and the growth rate of the functions that are provably total
in~$\SS$. For instance, the functions that are provably total in~$\IS1$ are the
primitive recursive functions, defined to be the functions which can be obtained
from the constants and the addition using the operations of composition and
definition by simple recursion.

The functions~$\Ack_\rr$ and $\Ack$ mentioned in
Section~\ref{S:Transition} are the functions defined by the following double
recursion rules: $\Ack_0(\xx) = \xx+1$, $\Ack_\rr(0) =
\Ack_{\rr-1}(1)$, and  $\Ack_\rr(\xx+1) =
\Ack_{\rr-1}(\Ack_\rr(\xx))$ for $\rr \ge 1$. For each~$\rr$, the
function~$\Ack_\rr$ is primitive recursive, hence provably total
in~$\IS1$, but the Ackermann function~$\Ack$ defined by~$\Ack(\xx)
= \Ack_\xx(\xx)$ eventually dominates all primitive recursive
functions and therefore is not provably total in $\IS1$. So, in order to
prove that a certain sentence~$\For$ is not provable from the axioms
of~$\IS1$, it is sufficient to establish that, from~$\For$, and using
arguments that can be formalized in~$\IS1$, one can prove the
existence of a function that grows as fast as the Ackermann function.

Let us mention that $\IS1$ is closely connected with the system~$\PRA$ of Primitive Recursive Arithmetic.
The latter is expressed in a language that contains the equality symbol and a symbol for each primitive 
recursive function, and its axioms essentially the defining equations for primitive recursive functions plus the
induction schema on formulas with bounded quantifiers. Primitive recursive functions are the provably total
functions of $\PRA$. As $\PRA$ and $\IS1$ turn out to prove the same $\Pi^0_2$ formulas, \ie, the same
aithmetical formulas of the form $\forall\xx_1 \exists\xx_2 (\For(\xx_1,\xx_2,\xx))$ with~$\For$ containing
bounded quantifiers only, they have the same provably total functions. 

As for provably total functions for the systems $\IS\kk$ with $\kk>1$
and for $\PA$, they can be characterized in terms of the Hardy Hierarchy of
fast-growing functions introduced in Definitions~\ref{D:Hardy3}
and~\ref{D:HardyN} using fundamental sequences of ordinals
below~$\eps_0$. This is a hierarchy of recursive functions $(H_\a)_{\a < \eps_0}$
indexed by transfinite ordinals and generalizing the Ackermann function, which
occurs as~$H_{\om^\om}$.  Then, a recursive function is
provably total in $\IS\kk$ (\resp in~$\PA$) if and only if it is primitive
recursive in\footnote{\ie, it can be obtained from the
constants, the projections, addition, and $H_\a$ using composition and definition by
simple recursion}~$H_{\alpha}$ for some
$\alpha <\om^{\om^{.^{.^\om}}}$,
$\kk+1$~times~$\om$ (\resp $\a < \eps_0$), a characterization that is 
extremely useful in proving independence  results from $\PA$ and its
subsystems~\cite{Min, Par}. 

\goodbreak
\begin{center}
\sc Acknowledgement
\end{center}

The authors warmly thank Andrey Bovykin who participated in early stages of this work but
could not join their group subsequently. They also thank Mireille Bousquet-Melou and Jean
Mairesse for useful suggestions.

\end{document}